\input  amstex
\input amsppt.sty
\magnification1200
\vsize=23.5truecm
\hsize=16.5truecm
\NoBlackBoxes

\def\supp{\operatorname{supp}}

\def\C{C^{\infty}}

\def\crp{\overline{\Bbb R}_+}

\def\rn{{\Bbb R}^n}
\def\rnp{{\Bbb R}^n_+}

\def\crnp{\overline{\Bbb R}^n_+}

\def\comega{\overline\Omega }

\def\ang#1{\langle {#1} \rangle}

\def\Pfrac{\tsize\frac1{\raise 1pt\hbox{$\scriptstyle p$}}}
\def\pfrac{\frac1{\raise 1pt\hbox{$\scriptscriptstyle p$}}}
\def\Pfracc#1{\tsize\frac{#1}{\raise 1pt\hbox{$\scriptstyle p$}}}
\def\pfracc#1{\frac{#1}{\raise 1pt\hbox{$\scriptscriptstyle p$}}}

\def\simto{\overset\sim\to\rightarrow}

\def\Zfrac{\tsize\frac1{\raise 1pt\hbox{$\scriptstyle z$}}}
\def\zfrac{\frac1{\raise 1pt\hbox{$\scriptscriptstyle z$}}}

\def\rp{ \Bbb R_+}

\def\rmi{ \Bbb R_-}

\def\OP{\operatorname{OP}}
\def\N{\Bbb N}
\def\R{\Bbb R}
\def\C{\Bbb C}
\def\Z{\Bbb Z}
\def\ol{\overline}

\def\F{\Cal F}

\document

\topmatter
\title The structure of solution spaces for fractional-order
operators, with gradient estimates \endtitle
\author Gerd Grubb \endauthor
\affil
{Department of Mathematical Sciences, Copenhagen University,
Universitetsparken 5, DK-2100 Copenhagen, Denmark.
E-mail {\tt grubb\@math.ku.dk}}\endaffil
\rightheadtext{Solution spaces}
\abstract
The solution space of the homogeneous Dirichlet problem for the
fractional Laplacian
$(-\Delta )^{a}$ ($0<a<1$) or a pseudodifferential generalization $P$,
on a bounded open set $\Omega \subset \rn$ with $C^{1+\tau }$-boundary,
$$
Pu=f \text{ on }\Omega ,\quad u=0 \text{ on }\rn\setminus \Omega ,
$$
is analysed in detail. It is shown, both for solutions in Sobolev
spaces of Bessel-potential type $H_q^t$ and  in H\"older-Zygmund spaces
$C_*^t$, that the solution space for $f$ of regularity $s\in [0,\tau -2a)$ is the
direct sum of a component $\dot H_q^{2a+s}(\comega)$ resp.\ $\dot
C_*^{2a+s}(\comega)$ with full regularity $2a+s$ and a component of
the form $d^a$ times a lifting of boundary values by Poisson
operators. Here $d(x)=\operatorname{dist}(x,\partial\Omega )$. This
extends to non-smooth problems results known in the $C^\infty $ setting. 

The knowledge is used to establish gradient estimates for $a>\frac12$,
such as e.g.,
$$
\aligned
\|d^{1-a+s}\nabla (u/d^a)\|_{\dot H_q^{s+t} (\comega )}&\le C
(\|f\|_{\ol H_q^t(\Omega )}+\|u\|_{L_q(\Omega )}),\\
\|d^{1-a+s}\nabla (u/d^a)\|_{\dot C_*^{s-\varepsilon } (\comega )}&\le C
(\|f\|_{L_\infty (\Omega )}+\|u\|_{L_\infty (\Omega )}),
\endaligned
$$
for small $s,t$,
with $\varepsilon =0$ if $\Omega $ and $P$ are $C^\infty $. The
estimates are entirely new in
Bessel-potential spaces; they extends previous results by Fall and
Jarohs in H\"older spaces.

A new tool is introduced: $\dot
 H^{s+t}_q(\comega)\subset d^s\dot H^{t}_q(\comega)$ holds for $s,t\ge
 0$ with $s+t<\tau $.

\endabstract
\endtopmatter

\subhead 1. Introduction \endsubhead

The regularity of solutions of 
 the
homogeneous Dirichlet problem with $f$ given in $H_q^s(\comega )$,
 $C^s(\comega )$ or $L_\infty (\Omega )$, $\Omega $ bounded open $ \subset \rn$,
$$
Pu=f\text{ on }\Omega , \quad u=0  \text{ on }\rn\setminus\Omega ,\tag1.1
$$
for fractional-order operators such as $P=(-\Delta )^a$ ($0<a<1$) or
pseudodifferential $x$-dependent generalizations $P=\OP(p(x,\xi ))$ of
order $2a$, has been intensively studied. An important step forward was when Ros-Oton and
Serra \cite{RS14} in 2014 showed for $P=(-\Delta )^a$, $\Omega $ being  $C^{1,1}$,
that 
$$
f\in L_\infty (\Omega ) \implies u\in d^aC^\alpha (\comega),\tag1.2
$$
(for small $\alpha $, later lifted to $\alpha $ up to $a$); here
$d(x)=\operatorname{dist}(x,\partial\Omega )$. The
author showed in \cite{G15} for general $P$ a complete characterization
in Bessel-potential spaces ($\Omega $ being $C^\infty $):
$$
f\in L_q (\Omega ) \iff u\in H_q^{a(2a)}(\comega),\tag1.3
$$
(valid also with $f\in \ol H_q^s(\Omega )$, $u\in
H_q^{a(2a+s)}(\comega)$, $s\ge 0$). Here the {\it $a$-transmission
space} $H_q^{a(2a+s)}(\comega)$ is defined by a formula $\Lambda _+ ^{(-a)}e^+\overline
  H_q^{a+s}(\Omega)$ (more details below); it satisfies
$$
H_q^{a(2a+s)}(\comega)\subset \dot H_q^{2a+s}(\comega)+d^ae^+\ol
H_q^{a+s}(\Omega ) ,\tag1.4
$$
(slightly modified if $s+a-\frac1q\in \Z$). $e^+ $ indicates extension
by zero.
There are many more contributions to the regularity question, e.g.\
\cite{RSV17}, and notably
 \cite{AR20} lifting (1.2) to
arbitrarily high orders.  We also described the
solvability of (1.1) in
the H\"older-Zygmund-scale $C_*^s(\comega)$, in \cite{G14}. Our results in terms of transmission spaces were extended to
nonsmooth situations in the joint work with Abels \cite{AG23}.

In \cite{G19} from 2019, the contribution from
$H_q^{a(2a+s)}(\comega)$ that actually enters in the second component $d^a\ol
H_q^{a+s}(\Omega )$ in (1.4) was characterized more precisely as $d^a$
times a space defined by a Poisson operator $\Cal K$ acting on weighted
boundary values of $u$, e.g.\ for low $s$,
$$
\aligned
u \in H_q^{a(2a+s)}(\comega)&\iff u=w+z,\; w\in \dot
H_q^{2a+s}(\comega),\; z\in d^aK_{(0)}B_q^{a+s-1/q}(\partial\Omega ),\\
u \in C_*^{a(2a+s)}(\comega)&\iff u=w+z,\; w\in \dot
C_*^{2a+s}(\comega),\; z\in d^aK_{(0)}C_*^{a+s}(\partial\Omega ),
\endaligned
\tag1.5
$$
for
$C^\infty $-domains $\Omega $.

Fall and Jarohs \cite{FJ21} showed in 2021 some interesting gradient estimates for
solutions $u$ in case $P=(-\Delta )^a$, $a\in (\frac12,1)$, e.g. 
$$
\|d^{1-a}\nabla u\|_{C^\beta (\comega)}\le C(\|f\|_{L_\infty (\Omega )}+\|u\|_{L_\infty (\Omega )}),\tag1.6
$$
with $\beta \in (0,2a-1)$, $\Omega $ being $C^{1,1}$. They point out
that (1.6) does not follow from (1.2), nor from the description as in
(1.4).
In the case of $C^\infty $-domains $\Omega $, we can show
(1.6) by use of (1.5). We shall here pursue the results also for nonsmooth domains, 
and as a new subject, we address the gradient estimates  in the
Bessel-potential scale $H_q^s$.

The present paper has several purposes:

1) One is to extend the
description of the transmission spaces
as in (1.5)
from \cite{G19} to nonsmooth domains
$\Omega $.

2) Another
is to show how the desired gradient estimates can be obtained from this
 precise description. The results pertain also to
$x$-dependent pseudodifferential operators $P$ with 
H\"older-dependence on $x$.

3) A third purpose is to present a new embedding
theorem for the $H_q^s$-scale of function spaces in terms of powers of
$d$, which serves to get good results for the component of $u$ in $\dot H^{2a+s}(\comega)$.

We now list some typical cases of the main results.

On the weighted trace- and Poisson-type-operators connected with
$a$-transmission spaces:

\proclaim{Theorem 1.1} Let $a\in (0,1)$, $\tau \ge 1$, $q \in
  (1,\infty )$, and  let $\Omega\subset\R^n$ be a bounded
  $C^{1+\tau}$-domain. Let $d(x)=\operatorname{dist}(x,\partial\Omega
  )$ near $\partial\Omega $, extended as a smooth function to
  the rest of $\Omega $.

The weighted trace operator $$
\gamma _0^au=\Gamma (a+1)\gamma_0(u/d^a)\colon H_q^{{a} (t)}(\comega )\to B_q^{t-a
  -1/q}(\partial \Omega ), \quad a+\tfrac1q <t<\tau +1,\tag1.7
$$
has a right inverse $K_{(0)}^a$ (such that $d^{-a}K^a_{(0)}$ is constructed in nonsmooth local coordinates from the
Poisson operator $K_0=\operatorname{OPK}((\ang{\xi '}+i\xi
_n)^{-1})$ on $\rnp$), satisfying
$$
K_{(0)}^{a}\colon B_{q}^{t-a
  -\frac1q}(\partial \Omega ) \to H_q^{{a} (t)}(\comega
  )\cap e^+d^a\ol H_q^{t-a }(\Omega  ).\tag1.8
$$

There is a similar result in the H\"older-Zygmund scale $C_*^s$. 

\endproclaim

On decomposition formulas for $a$-transmission spaces:

\proclaim{Theorem 1.2}  Let $a,q,\tau ,\Omega ,d $ be as in Theorem
{\rm 1.1}.

 $1^\circ$  For $s\in (-a+\frac1q,-a+1+\frac1q)$ with $s<\tau +1-2a$,
$H_q^{a(2a+s)}(\comega)$
is a direct sum:
$$
H_q^{a(2a+s)}(\comega)=\dot H_q^{2a+s}(\comega)\; \dot+ \;
  K^a_{(0)}B_q^{a+s-1/q}(\partial\Omega ).
\tag1.9
$$
With elements denoted $u=w+z$, $w$ runs through $\dot H_q^{2a+s}(\comega)$, and
$z$ equals $K^a_{(0)}\psi $, where
$\psi $ runs through $B^{a+s-1/q}_q(\partial\Omega )$. Moreover,
$\psi =\gamma ^a _0u$.

$2^\circ$ For $s\in (-a,1-a)$  with $s<\tau +1-2a$, $C_*^{a(2a+s)}(\comega)$
is a direct sum:
$$
C_*^{a(2a+s)}(\comega)=\dot C_*^{2a+s}(\comega)\; \dot+ \;
  K^a_{(0)}C_*^{a+s}(\partial\Omega ).
\tag1.10
$$
With elements denoted $u=w+z$, $w$ runs through $\dot C_*^{2a+s}(\comega)$, and
$z$ equals $K^a_{(0)}\psi $, where
$\psi $ runs through $C_*^{a+s}(\partial\Omega )$. Moreover,
$\psi =\gamma ^a _0u$. 
\endproclaim

On embeddings of
Bessel-potential spaces with powers of the
distance function:

\proclaim{Theorem 1.3}
 Let $s>0$,  $t\ge 0$. When $\Omega $ is a bounded $C^{1+\tau
 }$-domain with $\tau
\ge 1$, and  $t,s\ge 0$, one has for $s+t<\tau $ that $\dot
 H^{s+t}_q(\comega)\subset d^s\dot H^{t}_q(\comega)$; in other words: 
$$
\|d^{-s}u\|_{\dot H_q^{t}(\comega)}\le C \|u\|_{\dot H^{s+t}_q(\comega)}.\tag1.11
$$
Moreover, for $s,t<\tau $,
$$
\|d^{s}u\|_{\dot H_q^{t}(\comega)}\le C \|u\|_{\dot H^{t}_q(\comega)}.\tag1.12
$$
\endproclaim

This theorem extends  to the
 Bessel-potential scale some embedding properties known for H\"older spaces.
Notably, the first term in (1.4) can now be absorbed
 in the second term:

\proclaim{Corollary 1.4} Let $a,q,\tau ,\Omega ,d $ be as in Theorem
{\rm 1.1}. For $-a+\frac1q<s<\tau -2a$, the inclusion {\rm (1.4)} simplifies to
$$
H_q^{a(2a+s)}(\comega)\subset d^ae^+\ol
H_q^{a+s}(\Omega ) .\tag1.13
$$

\endproclaim

With the above
 tools available, we can show gradient estimates for solutions of
 (1.1), based on the knowledge of the solution spaces $H_q^{a(2a+s)}(\comega)$
and $C_*^{a(2a+s)}(\comega)$. Let $P$ be a pseudodifferential operator of order $2a$
 ($0<a<1$), classical, even and strongly elliptic, with symbol
 depending $C^\tau $ on $x$, $\tau >0$. Let $\Omega \subset \rn$ be  bounded
 open with $C^{1+\tau }$-boundary. The following results in Bessel-potential spaces are entirely new:

\proclaim{Theorem 1.5} Let $a\in (\frac12,1)$, $q\in (1,\infty )$,
$\tau \ge 1$, $t\in [0,\tau -2a)$.

When $u$ is a solution of {\rm (1.1)} with $f\in \ol
 H_q^t(\Omega )$, then
 $$
\|d^{1-a}\nabla u\|_{\ol H_q^{2a-1+t } (\Omega )^n}\le C (\|f\|_{\ol
H_q^t(\Omega )}+\|u\|_{L_q(\Omega )}).
\tag1.14
$$

If  moreover $s\in [0,2a-1)$ and $t\in [0,1-a+\frac1q)$ with
$s+t<\tau -2a$, then
$$
\|d^{1-a+s}\nabla (u/d^a)\|_{\dot H_q^{s+t} (\comega )}\le C
(\|f\|_{\ol H_q^t(\Omega )}+\|u\|_{L_q(\Omega )}).
\tag1.15
$$

When there is uniqueness of solution, the terms  $\|u\|_{L_q
(\Omega)}$ can be omitted.
\endproclaim

In H\"older-Zygmund spaces $C_*^s$, we have the following results,
extending and generalizing \cite{FJ21}:

\proclaim{Theorem 1.6} Let $a\in (\frac12,1)$, $s\in [0,2a-1)$, $t>0$,
and assume that $s+t<\tau -2a$. Then the solutions of {\rm (1.1)} satisfy
for small $\varepsilon $:
$$
\aligned
\|d^{1-a}\nabla u\|_{\ol C_*^{s+t-\varepsilon } (\Omega  )^n}&\le C
(\|f\|_{\ol C_*^t(\Omega )}+\|u\|_{L_\infty (\Omega )}),\\
\|d^{1-a}\nabla u\|_{\ol C_*^{2a-1-\varepsilon } (\Omega  )^n}&\le C
(\|f\|_{L_\infty (\Omega )}+\|u\|_{L_\infty (\Omega )}).
\endaligned
\tag1.16
$$

When moreover  $t\in (0,1-a)$, $s>0$,
$$
\aligned
\|d^{1-a+s}\nabla (u/d^a)\|_{\dot C_*^{s+t-\varepsilon } (\comega )}&\le C
(\|f\|_{\ol C_*^t(\Omega )}+\|u\|_{L_\infty (\Omega )}),\\
\|d^{1-a+s}\nabla (u/d^a)\|_{\dot C_*^{s-\varepsilon } (\comega )}&\le C
(\|f\|_{L_\infty (\Omega )}+\|u\|_{L_\infty (\Omega )}).
\endaligned
\tag1.17
$$

The results hold  with ``$-\varepsilon
$'' removed when $\Omega $ and $P$ are $C^\infty $ in $x$.
\endproclaim

\subsubhead Plan of the paper \endsubsubhead
The function spaces and operators are introduced in Section 2. In
Section 3, we recall the decomposition results for solution spaces known in the smooth
situation. Section 4 begins the new material, establishing the precise
description of solution spaces in nonsmooth situations. Section 5
shows a useful new result on embeddings in Bessel-potential spaces
(and other function scales) involving powers of the distance to the
boundary.
In Section 6 it is shown how the preceding analysis allows to show
gradient estimates for solutions of the homogeneous Dirichlet problem,
both in Bessel-potential spaces and H\"older-Zygmund spaces. The
Appendix shows a needed multiplication theorem for $\mu $-transmission spaces
over nonsmooth sets.

\subsubhead Acknowledgement \endsubsubhead The author is grateful to
 Helmut Abels and Hans Triebel for useful conversations.

\subhead 2. Preliminaries \endsubhead

\subsubhead 2.1 Function spaces \endsubsubhead

In the following, $\Omega \subset \rn$ is an open set with $C^\infty $
or $C^{1+\tau }$-boundary, $\tau >0$. We assume throughout the paper that $q\in (1,\infty )$.

Recall that the Sobolev spaces of Bessel-potential type
$H^s_q$ are defined, with  $\ang\xi =(1+|\xi
|^2)^{\frac12}$, $s\in\Bbb R$, by
$$
\aligned
H_q^s(\rn)&=\{u\in \Cal S'(\Bbb R^n)\mid \Cal F^{-1}(\langle{\xi }\rangle^s\hat
u)\in L _q(\Bbb R^n)\},\\
  \overline H_q ^s(\Omega  )&=r^+H_q^s (\Bbb R^n),\text{ the
  restricted space},\\
\dot H_q^s (\comega )&=\{u\in H_q^s (\Bbb R^n)\mid
\operatorname{supp}u\subset \comega  \},\text{ the supported space}.
\endaligned \tag2.1
$$
Here $r^+$ denotes restriction to $\Omega $; $e^+$ will indicate
extension by 0 from $\Omega $ to $\rn$.
$\overline H_q^s(\Omega )$ and $\dot H_{q'}^{-s}(\comega)$ are dual
spaces, with a duality consistent with the $L_2$-scalar
product ($\frac1q+\frac1{q'}=1$).

For $q=2$,  the index can be omitted. The dot and overline notation
 stems from H\"ormander's books.

The spaces are in some texts denoted $L^{s,q}$ or $L^s_q$ and defined only
for $s\ge 0$ (following a notation in Stein's book \cite{S70}); the notation we use
here is consistent with many other works and in particular Triebel's books on function spaces and
interpolation theory, e.g.\ \cite{T95}. They equal the
$W^{s,q}$-spaces when $s\in\N_0$.

Besov spaces will not play an important role in this paper, except that
they enter naturally as spaces of boundary values; we recall that the
mapping $\gamma _0\colon u(x)\mapsto u(x',0)$ is continuous and
surjective from $\ol H_q^t(\rnp)$ to $B_{q}^{t-1/q}(\R^{n-1})$ for
$t>\frac1q$ (here $B_q^{s}(\R^{n-1})$ is also known as
$B_{qq}^{s}(\R^{n-1})$; it equals
 $W^{s,q}(\R^{n-1})$ when $s\in \rp\setminus \N$).

The H\"older-Zygmund spaces $C_*^t(\rn)$ coincide with the H\"older spaces $C^t(\rn)$ for $t\in\rp \setminus
\N$; they are also defined for negative values $t$ and have nice
interpolation properties on the whole scale. They identify with the
spaces $B^s_{\infty ,\infty }(\rn)$ in the Besov-Triebel-Lizorkin
scales of function spaces (more details e.g.\ in  \cite{G14}).  There are spaces over
$\Omega $:
$$
\aligned
  \overline C_* ^s(\Omega  )&=r^+C_*^s (\Bbb R^n),\text{ the
  restricted space},\\
\dot C_*^s (\comega )&=\{u\in C_*^s (\Bbb R^n)\mid
\operatorname{supp}u\subset \comega  \},\text{ the supported space}.
\endaligned \tag2.2
$$
In particular, $C^0(\comega)\subset L_\infty (\Omega )\subset \ol
C^0_*(\Omega)$ when $\Omega $ is bounded.

\subsubhead 2.2 The $\mu $-transmission spaces \endsubsubhead

We consider a pseudodifferential generalization $P$ of the fractional
Laplacian $(-\Delta )^a$ (including $(-\Delta )^a$ as a special case),
defined from a smooth symbol as in \cite{G15} or a nonsmooth symbol as
in Abels-Grubb \cite{AG23}, $0<a<1$. (The power is denoted $a$ where
some other works use $s$; this is because
$s$ is heavily used as a differentiability index on function spaces.)

Recall that a {\it pseudodifferential operator} 
$P$ on ${\Bbb R}^n$ is
defined from a symbol $p(x,\xi )$ on ${\Bbb R}^n\times{\Bbb R}^n$ by 
$$
Pu=\operatorname{OP}(p(x,\xi ))u 
=(2\pi )^{-n}\int _{\rn}e^{ix\cdot\xi
}p(x,\xi )\hat u(\xi )\, d\xi =\Cal F^{-1}_{\xi \to x}(p(x,\xi )\F u(\xi
)),\tag 2.3
$$  
using the Fourier transform $\F\colon u(x)\mapsto \hat u(\xi )=\int
e^{-ix\cdot \xi }u(x)\, dx$. The formula makes sense for nice
functions $u$, and is generalized to distributions by interpretations
as ``oscillatory integrals'' (references to the literature are given
e.g.\ in \cite{G19}, \cite{AG23}).
Here the symbols $p(x,\xi )$ are
generally taken to lie in the symbol space $S^m_{1,0}({\Bbb R}^n\times{\Bbb R}^n)$, consisting of
$C^\infty $-functions $p(x,\xi )$
such that $\partial_x^\beta \partial_\xi ^\alpha p(x,\xi
)$ is $O(\ang\xi ^{m-|\alpha |})$ for all $\alpha ,\beta $, for some
$m\in{\Bbb R}$, with global estimates for $x\in{\Bbb R}^n$.  $P$ (of order $m$)
then  maps $H^s_q({\Bbb R}^n)$ continuously into
$H^{s-m}_q ({\Bbb R}^n)$ for all $s\in{\Bbb R}$. One adds to these
operators the ``negligible'' operators, that is, linear operators from
distributions to $C^\infty $; they occur in the theory as error terms
(remainders) in any case.

We shall also consider symbols that are only in the H\"older class
$C^\tau $ with respect to $x$ for some $\tau \in (0,\infty )$; the
space
$C^\tau S^m_{1,0}(\rn\times\rn)$ consists of functions $p(x,\xi )$
satisfying estimates 
$$
\|\partial_\xi ^\alpha p(\cdot,\xi )\|_{C^\tau (\rn)}\le C_\alpha \ang{\xi
}^{m-|\alpha |}, \text{ all }\alpha \in \N_0^n.         
$$
Here  $S^m_{1,0}({\Bbb R}^n\times{\Bbb R}^n)$ is included as
$C^\infty S^m_{1,0}({\Bbb R}^n\times{\Bbb R}^n)$; the case $\tau =\infty $.

\proclaim{Hypothesis 2.1} $P$ is a  ps.d.o.\
$P=\operatorname{OP}(p(x,\xi ))$ of order $2a$, $a\in(0,1)$, where $p$
belongs to 
$C^\tau
S^{2a}_{1,0}(\rn\times\rn)$ for some $\tau \in (0,\infty ]$, and satisfies
(with formulas  valid for $|\xi
|\ge 1$):

{\rm (i)} $p$ is {\bf
  classical:}
$p\sim \sum_{j\in
  {\N}_0}p_j$ with $p_j(x,t\xi )=t^{2a-j}p_j(x,\xi )$;

{\rm (ii)} $p$ is {\bf even:}  $p_j(x,-\xi )=(-1)^jp_j(x,\xi )$, all
$j$;

{\rm (iii)}
$p$ is {\bf strongly elliptic:} $\operatorname{Re}p_0(x,\xi
)\ge c|\xi |^{2a}$, with $c>0$. 

\endproclaim

The full pseudodifferential calculus includes smoothing operators. The
fractional Laplacian $(-\Delta )^a$ fits into the picture in the way
that  $(-\Delta )^a=\OP(|\xi |^{2a})=P_1+P_2$, where
$P_1=\OP((1-\psi(\xi ))|\xi |^{2a} )$, with $\psi \in C_0^\infty(\rn)$
equal to 1 for $|\xi |\le \frac12$, is as in Hypothesis 2.1; and
$P_2=\OP(\psi(\xi )|\xi |^{2a} )$ is smoothing,  mapping $\bigcup_{t\in\R}
H_q^t(\rn)$, hence $\Cal E'(\rn)$, to $\bigcap_{t\in\R}
H_q^t(\rn)\subset C^\infty (\rn)$. Using this,  $P=(-\Delta )^a$ is
included in the regularity analysis over bounded
sets. Also fractional powers of strongly elliptic second-order
differential operators with smooth coefficients are included.

In (i), the asymptotic property $p\sim \sum_{j\in
  {\N}_0}p_j$ indicates that for all $J\in \N$, $p-\sum_{j<J}p_j$ is
in $S^{2a-J}_{1,0}(\rn\times\rn)$ resp.\ $C^\tau
S^{2a-J}_{1,0}(\rn\times\rn)$. The spaces of {\it classical} symbols are
  denoted $S^{2a}(\rn\times\rn)$ resp.\ $C^\tau
S^{2a}(\rn\times\rn)$ without the $\{1,0\}$ subscript.

Condition (ii) means that $p$ has the {\it $a$-transmission
  property} (from \cite{G15}, \cite{AG23}) for all considered domains $\Omega $.

For such operators $P$, the {\bf homogeneous Dirichlet problem} over an
open 
subset  $\Omega \subset \rn$
is: Given a function $f$ on $\Omega $, find $u$ such that
$$
Pu= f \text{ in }\Omega ,\quad u=0 \text{ on }\rn\setminus \Omega .\tag2.4
$$ 
For clarity, one can here write ``$r^+Pu=f$''.
It was shown for smooth bounded domains $\Omega $ and smooth $P$ in \cite{G15}, Th.\ 4.4, that (taking $u\in
\dot H^{a-\delta } _q(\comega )$ for some $\delta <\frac1{q'}$) there
holds when $s>-a-\tfrac1{q'}$:
$$
f\in \overline H^s_q(\Omega )\iff u\in H^{a(2a+s)}_q
(\overline\Omega ).\tag2.5
$$
\cite{AG23} extended this to  finite $\tau >2a$,  $s\in [0,\tau -2a)$,
when $\Omega $ is $C^{\tau +1}$ and $P$ has symbol $p(x,\xi )$ in $C^\tau
S^{2a}(\rn\times\rn)$.

Here the spaces $ H^{a(2a+s)}_q
(\overline\Omega )$   are so-called
{\bf $a$-transmission spaces}, defined in the $L_2$-case in an unpublished note by
H\"ormander from 1966 \cite{H66}, and introduced more generally in \cite{G15} for
smooth domains:
$$
H_q^{a(2a+s)}(\comega)=\Lambda _+ ^{(-a)}e^+\overline
  H_q^{a+s}(\Omega),\text{ for }s>-a-\tfrac1{q'}.\tag2.6
$$
More about them below; we shall first formulate solvability statements
resulting from \cite{AG23,G23}.

The problem (2.4) gives rise to the definition of the Dirichlet
realization $P_{D,q}$ in $L_q(\Omega )$, acting as $r^+P$ and with domain
$$
D(P_{D,q})=\{u\in \dot H_q^a(\comega)\mid r^+Pu\in L_q(\Omega )\},
\text{ equal to }H_q^{a(2a)}(\comega)\tag2.7
$$
in view of (2.5) with $s=0$. It is accounted for e.g.\ in \cite{G23, Sect.\ 4} that
$P_{D,2}$ in $L_2(\Omega )$ is lower bounded with compact resolvent,
having its discrete spectrum $\Sigma $ contained in a sectorial region
$M$:
$$
\Sigma \subset M=\{\lambda \in\C \mid \operatorname{Re}\lambda \ge c-\beta ,
|\operatorname{Im}\lambda |\le c_1(\operatorname{Re}(\lambda +\beta
)\}.\tag2.8
$$
Moreover, $P_{D,q}$ has for each $q\in (1,\infty )$ the same
spectrum $\Sigma $, with the same eigenspaces $N_\lambda $ (and the same
eigenspaces $N'_{\bar\lambda }$ of $(P_{D,q})^*=P^*_{D,q'}$ with $\operatorname{dim}N_\lambda =\operatorname{dim}N'_{\bar\lambda }<\infty $) as in
the case $q=2$; they
are contained in $\dot C^a(\comega)$. We now list some basic known facts:

\proclaim{Theorem 2.2} Let $q\in (1,\infty )$, $a\in (0,1)$, $\tau
>2a$, $s\in [0,\tau -2a)$. Let $\Omega $ be a bounded $C^{1+\tau
}$-domain in $\rn$, and let $P$ satisfy Hypothesis {\rm 2.1}.

If
$0$ is not in the spectrum $\Sigma $ of $P_{D,q}$, $r^+P$ defines a homeomorphism
$$
r^+P\colon H_q^{a(2a+s)}(\comega) \simto \ol H_q^s(\Omega);\tag2.9
$$
  hence for $u$ and $f$ in {\rm (2.4)},
$$
\|u\|_{H_q^{a(2a+s)}(\comega)}\simeq \|f\|_{\ol H_q^s(\Omega)}.\tag2.10
$$

Generally, there holds with positive constants $C,C'$:
$$
\|f\|_{\ol H_q^s(\Omega)}\le C\|u\|_{H_q^{a(2a+s)}(\comega)}\le C'
(\|f\|_{\ol H_q^s(\Omega)} +\|u\|_{L_q(\Omega)}).\tag2.11
$$
\endproclaim

\demo{Proof} (2.5) was extended to the nonsmooth situation in \cite{AG23} and
recalled in \cite{G23, Th.\ 3.2}; it clearly implies (2.9) and (2.10)
when  $0\notin \Sigma $.

Whether $0\in\Sigma $ or not, the left-hand inequality in (2.11) holds by the
forward mapping property of $r^+P$ recalled in \cite{G23, Th.\
3.2 $1^\circ$}. For the right-hand inequality, we have from \cite{G23, Prop.\
3.3} that
$$
\|u\|_{H_q^{a(2a+s)}(\comega)}\le C_1
\|f\|_{\ol H_q^s(\Omega)} +C_1\|u\|_{\dot H_q^a(\comega)}.\tag2.12
$$
Since for small $\delta >0$,
$$
\|u\|_{\dot H_q^a(\comega)}\le \varepsilon \|u\|_{\dot
H_q^{a+1/q-\delta }(\comega)} +C(\varepsilon ) \|u\|_{L_q(\Omega)}
\le C_0\varepsilon  \|u\|_{ H_q^{a(2a+s) }(\comega)}+C(\varepsilon )\|u\|_{L_q(\Omega)}
$$
with a fixed $C_0$ (cf.\ (2.17) below),
we can change $C_1\|u\|_{\dot H_q^a(\comega)}$ to
$C_2\|u\|_{L_q(\Omega)}$ in (2.12).\qed

\enddemo

We now recall the systematic definition of $\mu $-transmission spaces
in the smooth case where $\Omega $ is
either $\rnp$ or an open bounded $C^\infty $ subset of $\rn$.
Let $\mu >-1$\footnote{This will  primarily be used
with $\mu =a$ or $\mu =a-1$.} and let $t>-\mu -\frac1{q'}$. Then the
$\mu $-transmission spaces are defined by
$$
H_q^{\mu(t)}(\comega)=\Lambda _+ ^{(-\mu )}e^+\overline
  H_q^{t-\mu}(\Omega),\text{ for }t>-\mu-\tfrac1{q'},\tag2.13
$$
where $$
  \Lambda _+^{(m)}\colon \dot H_q^{r}(\comega) \simto
  \dot H_q^{r-m}(\comega),\quad r,m\in\R,\tag2.14
  $$
  is a family of ``order-reducing'' bijective pseudodifferential
  operators of order $m$ preserving support in
  $\comega$ (\cite{G90},\cite{G15}). 
  For $\Omega =\rnp$, one can use
  $
  \Xi _+^m=\OP((\ang{\xi
  '}+i\xi _n)^m) 
$ instead of $\Lambda _+^{(m)}$.

Now some properties established in \cite{G15}, \cite{G14} will be
listed: 
$$
H_q^{\mu(t)}(\comega)  =
 \dot
H_q^{t} (\overline\Omega ),\text{ when
 }-\tfrac1{q'}<t-\mu<\tfrac1q,\text{ in particular }H_q^{\mu(\mu )}(\comega)  =
 \dot
H_q^{\mu } (\overline\Omega ).\tag2.15
$$
When
  $t-\mu \ge \frac1q$, $\varepsilon >0$,
$$
\aligned
 H_q^{\mu(t)}(\comega)&\supset
 \cases  \dot H_q^{t} (\overline\Omega )\text{ when }t-\mu-\frac1q\notin
 \N_0,\\
\dot H_q^{t+\varepsilon } (\overline\Omega )\text{ when }t-\mu-\frac1q\in
 \N_0,\endcases\\
   H_q^{\mu(t)}(\comega)&\subset  \cases
 \dot
H_q^{t} (\overline\Omega )+e^+ d^\mu\overline H_q^{t-\mu}
                        (\Omega )\text{ when }t-\mu-\frac1q\notin \N_0,\\
  \dot
H_q^{t-\varepsilon } (\overline\Omega )+e^+ d^\mu\overline H_q^{t-\mu}
                        (\Omega )\text{ when }t-\mu-\frac1q\in \N_0.
                        \endcases
\endaligned \tag2.16
$$
When $\mu \ge 0$, the second statement in
(2.16) can thanks to Theorem 1.3 (Theorem 5.4 below) be simplified for $t>\mu +\frac1q$ to
$$
 H_q^{\mu(t)}(\comega)\subset e^+ d^\mu\overline H_q^{t-\mu}
                        (\Omega ), \tag2.17
$$
since
$
 \dot
H_q^{t} (\overline\Omega )\subset e^+ d^\mu\dot H_q^{t-\mu}$ (the
cases $t-\mu -\frac1q\in\N$ are included by interpolation).
 
For the case $t=2\mu \ge 0$:
$$
H_q^{\mu(2\mu)}(\comega)  \subset
 \dot
H_q^{\mu+\frac1q-\delta } (\overline\Omega ),\text{ any }\delta  >0.\tag2.18
$$

If $\Omega $ is only $C^{1+\tau }$, the spaces are considered
for  $t\in (-\frac1{q'},\tau -2\mu)$; here
$H_q^{\mu(t)}(\comega)$ is defined by a localized version of formula
(2.13), cf.\ \cite{AG23, Sect.\ 4.2}. More on this further below in
Section 4.

The set-up of $\mu $-transmission spaces in the H\"older-Zygmund scale is
completely analogous, cf.\ \cite{G14}.
The auxiliary ps.d.o.s $\Lambda _+^{(m)}$ map in a similar way as in
(2.14):
$$
  \Lambda _+^{(m)}\colon \dot C_*^{r}(\comega) \simto
  \dot C_*^{r-m}(\comega),\quad r,m\in\R,\tag2.19
  $$
In comparison with $H^t_q$-spaces, the rule of thumb is that $q$
is replaced by $\infty $, whereby 
$\frac 1q$ is replaced by $0$, and $\frac 1{q'}$ and $q'$ are replaced
by $1$.
In this scale, $\ol C_*^s(\Omega )$ identifies with $\dot C_*^s(\comega)$ for
$s\in (-1,0)$.

Let $\mu >-1$ and let $t>-\mu -1$. Then the $\mu $-transmission spaces
are defined by:
$$
C_*^{\mu(t)}(\comega)=\Lambda _+ ^{(-\mu )}e^+\overline
  C_*^{t-\mu}(\Omega),\text{ for }t>-\mu-1;\tag2.20
$$
 $\Xi ^{-\mu }_+$ can be used instead of $\Lambda _+ ^{(-\mu )}$ when $\Omega =\rnp$.

There are the identifications
and inclusions:
$$
C_*^{\mu(t)}(\comega)  =
 \dot
C_*^{t} (\overline\Omega ),\text{ for }-1<t-\mu<0,\tag2.21
$$
and when
  $t-\mu \ge 0$, $\varepsilon >0$,
$$
\aligned
 C_*^{\mu(t)}(\comega)&\supset
 \cases  \dot C_*^{t} (\overline\Omega )\text{ when }t-\mu\notin
 \N_0,\\
\dot C_*^{t+\varepsilon } (\overline\Omega )\text{ when }t-\mu\in
 \N_0,\endcases\\
   C_*^{\mu(t)}(\comega)&\subset  \cases
 \dot
C_*^{t} (\overline\Omega )+e^+ d^\mu\overline C_*^{t-\mu}
                        (\Omega )\text{ when }t-\mu\notin \N_0,\\
  \dot
C_*^{t-\varepsilon } (\overline\Omega )+e^+ d^\mu\overline C_*^{t-\mu}
                        (\Omega )\text{ when }t-\mu\in \N_0.
                        \endcases
                        \endaligned
\tag2.22
$$
The inclusion $\dot C_*^{s+t}(\comega)\subset d^s\dot
C_*^{t}(\comega)$ for $s,t\ge 0$ (cf.\ Section 5 below) allows to simplify the second statement in (2.22)
for $t>\mu \ge 0$ to
$$
  C_*^{\mu(t)}(\comega)\subset e^+ d^\mu\overline C_*^{t-\mu}
                        (\Omega ).\tag 2.23
$$
This fact was already used in \cite{G14} in particular cases, and
complies well with the findings in works of Ros-Oton and coauthors.

For pseudodifferential operators $P$, the extension of  mapping
properties and regularity properties from $C^\infty $-cases to $C^{\tau}$-cases
was in \cite{AG23} focused on the Bessel-potential scale
$H^s_q$. Results for other scales of function spaces can to some
extent be derived from those results, but a thorough investigation of
other scales in their own right has not been worked out in non-smooth cases, as far as we
know (it would be quite demanding). One can to some
extent derive results in H\"older-Zygmund spaces from results in
Bessel-potential spaces by approximation.

In the smooth case, results in H\"older-Zygmund spaces can be set up
directly using \cite{G14}.

Our systematic results for the problem (2.4) in H\"older-Zygmund
spaces are as follows, when $\mu =a\in (0,1)$:

\proclaim{Theorem 2.3} $1^\circ$ Let $a\in (0,1)$ and $s>a-1$, let $\Omega $ be a bounded $C^\infty
$-domain, and let $P$ be a ps.d.o.\ satisfying Hypothesis {\rm 2.1}
with $\tau =\infty $.  Then
$r^+P$ maps continuously
$C_*^{a(2a+s)}(\comega)$ to $\ol C_*^s(\Omega )$.
Moreover, if
$u\in \dot C_*^\sigma(\comega) $ (some $\sigma >a-1$) solves {\rm (2.4)} with $f\in
\ol C_*^{s}(\Omega )$, then $u\in C_*^{a(2a+s)}(\comega)$.
The mapping $r^+P$
defines a Fredholm
operator with index $0$ from $C_*^{a(2a+s)}(\comega)$ to $\ol C_*^s(\Omega )$. In the
case with uniqueness, $$
\|u\|_{C_*^{a(2a+s)}(\comega)}\simeq \|f\|_{\ol C_*^s(\Omega)};\tag2.24
$$
in general there holds:
$$
\|f\|_{\ol C_*^s(\Omega)}\le C\|u\|_{C_*^{a(2a+s)}(\comega)}\le C'
(\|f\|_{\ol C_*^s(\Omega)} +\|u\|_{L_\infty (\Omega)})\tag2.25
$$
(the last norm can be replaced by $\|u\|_{\dot C_*^\sigma(\comega)}$).

$2^\circ$ Let $a\in (0,1)$,  $\tau \in (0,\infty )$,  $s\in (0,\tau -2a)$;  let $\Omega $ be a bounded $C^{1+\tau }
$-domain, and let $P$ be a ps.d.o.\ satisfying Hypothesis {\rm 2.1}
for the given $\tau $. 
Let $u\in \dot H_q^a(\comega)$ for some $q\in (1,\infty )$.

If $u$ solves {\rm (2.4)} for an $f\in \ol C_*^{s}(\Omega )$,
then $u\in C_*^{a(2a+s-\varepsilon )}(\comega)$ for every small
$\varepsilon >0$.

If $u$ solves {\rm (2.4)} for an $f\in L_\infty (\Omega )$,
then $u\in C_*^{a(2a-\varepsilon )}(\comega)$ for every small
$\varepsilon >0$.
\endproclaim

\demo{Proof} The result in $1^\circ$ is a special case of Theorem 3.2
and Remark 3.3 in \cite{G14}. The estimates (2.25) follow as in
Theorem 2.1.

The result in $2^\circ$ is shown in \cite{AG23, Cor.\ 6.11}.
It is
derived from the result in $H_q^t$-spaces by use of Sobolev embedding,
noting that (with $t=2a+s$)
$$
H_q^{a(t)}(\comega) \subset C_*^{a(t-n/q-\varepsilon
)}(\comega),\text{ when }a\le t-n/q-\varepsilon <1+\tau , \tag2.26
$$
and $\ol C_*^s(\Omega )\subset \ol H_q^{s-\varepsilon /2}(\Omega )$.
\qed
 
\enddemo

\subhead 3. Decomposition results known for $C^\infty $-domains \endsubhead

\subsubhead 3.1  Results in Bessel-potential spaces  \endsubsubhead

In the following, we describe a more detailed analysis of the
$\mu $-transmission spaces, aimed at the cases $\mu =a$ and $\mu =a-1$
($0<a<1$). It was worked out for $C^\infty
$-domains $\Omega $ in \cite{G19} (preceded by observations in
\cite{G15, G18}). We here recall the $C^\infty $-results; they will be extended to
$C^{1+\tau }$-domains in Section 4. 

An important point is to show which elements of $e^+ d^\mu \overline H_q^{t-\mu }(\Omega )$ resp.\
$e^+ d^\mu \overline C_*^{t-\mu }(\Omega )$ in the decomposition
formulas (2.16) and (2.22) that actually enter in the space.
The results for $H^t_q$-spaces will be described first.

Let $\Omega \subset \rn$ be bounded and open with $C^\infty
$-boundary, or $\Omega =\rnp$. For systematic purposes, it is natural to provide the weighted trace maps
$\gamma _j^\mu $ with a gamma-function factor:
$$
\gamma _j^\mu u=\Gamma (\mu +1+j)\gamma _j(u/d^\mu );\tag3.1
$$
here $\gamma _0v=v|_{\partial\Omega }$ and  $\gamma _jv=(\frac{\partial}{\partial \nu
})^jv|_{\partial\Omega }$ for $j\in\N$; $d(x)=\operatorname{dist}(x,\partial\Omega
)$).

Along with trace operators $\gamma _j$ from $\comega $ to
$\partial\Omega $ enter also Poisson operators going from
$\partial\Omega $ to $\comega$, defined in the Boutet de Monvel
calculus (\cite{B71}, $L_p$-theory in \cite{G90}, Johnsen \cite{J96}
and later works). They are used in the present paper on an elementary level.
We first recall the construction of the auxiliary Poisson operator
$K_{(0)}$ of order 0, and its relation to weighted traces.  The basic calculation takes place in the model situation on
$\rnp$:

\proclaim{Lemma 3.1} 
Let $K_0$ denote the Poisson operator from ${\Bbb
R}^{n-1}$ to $\crnp$ with symbol $\frac 1{\ang{\xi '}+i\xi _n}$,
$$
K_0\varphi =\operatorname{OPK} (\tfrac1{\ang{\xi '}+i\xi _n})\varphi =\F^{-1}_{\xi
\to x}\tfrac {\hat \varphi (\xi ')}{\ang{\xi '}+i\xi _n}.
$$

When   $t>\mu +\tfrac1q $, the elements of $H_q^{\mu (t)}(\crnp)$ have a
unique decomposition
$$
\aligned
u&=w+z ,\text{ where }w \in H_q^{(\mu +1)(t)}(\crnp),\text{ and}\\
z &=e^+\tfrac1{\Gamma (\mu +1)}x_n^{\mu } K_0 \gamma ^\mu _0u\in e^+x_n^\mu \ol H_q^{t-\mu }(\rnp )\cap H_q^{\mu (t)}(\crnp).
\endaligned
$$
In other words, the elements of $H_q^{\mu (t)}(\crnp)$ are parametrized as 
$$
u=w+ K_0^\mu \varphi ,\quad K_0^\mu =  e^+\tfrac1{\Gamma (\mu +1)}x_n^{\mu } K_0 ,
$$
where 
$w $ runs through $ H_q^{(\mu +1)(t)}(\crnp)$ and $\varphi  $ runs through $B_q^{t-\mu -1/q}({\Bbb R}^{n-1} )$; here $\varphi  $ equals $\gamma ^\mu _0u$.

\endproclaim 

\demo{Remarks on the proof}
There is a detailed proof in \cite{G19, Lemma 3.3}. As a supplement to
the proof, we note: There holds that $e^+x_n^\mu
\Cal S(\crnp)$ is dense in $H_q^{\mu (t)}(\crnp)$, which
follows straightforwardly from the fact shown in \cite{G21, Lemma 6.1}: $e^+x_n^\mu
\Cal S(\crnp)= \Xi _+^{-\mu }e^+\Cal S(\crnp)$.

It should be noted that the indication $e^+$ in the formula for $K^\mu
_0$ is in a sense superfluous, since $K_0$ maps into a space of
distributions on $\rn$ supported in $\crnp$; it may be written there
for explicitness. 
\enddemo

The case of bounded $\Omega $ is treated in \cite{G19}, Th.\ 3.4:

\proclaim{Theorem 3.2}  Let $\Omega \subset {\Bbb R}^n$, bounded, open
with $C^\infty $-boundary. 

There is a
Poisson operator $ K_{(0)}$ of order $0$
from $\partial\Omega $ to $\comega$ 
(in the Boutet de Monvel calculus) with principal symbol $1/(|\xi
'|+i\xi _n)$ in local coordinates at the boundary,
such that $K_{(0)}$ is a right inverse of $\gamma _0$, and the
following holds:

The operator $K^\mu  _{(0)}$ defined by
$$
 K^\mu  _{(0)}= \tfrac1{\Gamma (\mu +1)}d^\mu  e^+ K_{(0)},\tag3.2
$$
 maps  $B_q^{t-\mu -1/q}(\partial\Omega )$ into $ e^+d^\mu  \ol
H_q^{t-\mu }(\Omega )\cap H_q^{\mu (t)}(\comega)$ for $t>\mu +\tfrac1q $,
and satisfies
$$
\gamma ^\mu  _0 K^\mu  _{(0)}\varphi =\varphi,\text{ all }\varphi \in B_q^{t-\mu -1/q}(\partial\Omega ) .\tag 3.3
$$

When $t>\mu +\tfrac1q $, the elements $u$ of $H_q^{\mu (t)}(\comega)$ have a
unique decomposition
$$
u=w+K^\mu _{(0)} \varphi ,\tag 3.4
$$
where 
$w  $ runs through $ H_q^{(\mu +1)(t)}(\comega)$ (equal to $\dot
H_q^t(\comega)$ when $t-\mu \in(\tfrac1q ,1+\tfrac1q )$) and $\varphi $ runs through $B_q^{t-\mu -1/q}(\partial\Omega )$; here  $\gamma ^\mu 
_0u$ equals $\varphi $.

\endproclaim

In the rest of this section,
 $\Omega \subset \rn $ is either bounded, open and $C^\infty$, or it equals $\rnp$. We consider in particular
the cases $\mu =a$ or $a-1$, $a\in (0,1)$. $K_{(0)}$ is to be read as
$K_0$ when $\Omega =\rnp$.

In the case $\mu =a$, the general result has the following simple corollary: 
Recall from \cite{G15} that the solution space for the homogeneous Dirichlet problem (2.4) with $f\in \ol H_q^s(\Omega )$
equals the $a$-transmission space $H_q^{a(2a+s)}(\comega)$. ($2a+s$
plays the same role as $t$ above.) For this space we have in particular
when 
$a+s\in (\frac1q,1+\frac1q)$ (as noted in \cite{G19, Th.\ 3.4}):

\proclaim{Corollary 3.3} Let 
$s\in
(-a+\frac1q, -a+1+\frac1q)$. The $a$-transmission
space $H_q^{a(2a+s)}(\comega)$  
is a direct sum:
$$
H_q^{a(2a+s)}(\comega)=\dot H_q^{2a+s}(\comega)\; \dot+ \;
  d^aK_{(0)}B_q^{a+s-1/q}(\partial\Omega ).
\tag3.5
$$
Denote the elements $u=w+z$, then
$w$ runs through $\dot H_q^{2a+s}(\comega)$, and
$z$ equals $d^aK_{(0)}\psi $, with
$\psi $ running through $B^{a+s-1/q}_q(\partial\Omega )$. Moreover, $\psi =\gamma _0(u/d^a)$.
\endproclaim

With the normalization (3.1) and (3.2), the formula for $z$ is
 $z=K^a_{(0)}\psi' $, $\psi' =\gamma ^a_0u$.

 Note that $w$ has the best possible regularity $2a+s$, whereas $z$
  describes exactly where the singularity  $d^a$ comes from.

For nonhomogeneous boundary value problems, we must go to the
transmission spaces with $\mu =a-1$, cf.\ e.g.\ \cite{G14},\cite{G18}. Recall
from (2.6), (2.15), (2.16) that 
$$
H_q^{(a-1)(t)}(\comega)=\Lambda _+ ^{(-a+1 )}e^+\overline
  H_q^{t-a+1}(\Omega),\text{ for }t>-a+1-\tfrac1{q'},\tag3.6
$$
satisfying
$$
H_q^{(a-1)(t)}(\comega)  =
 \dot
H_q^{t} (\overline\Omega ),\text{ for }-\tfrac1{q'}<t-a+1<\tfrac1q,\tag3.7
$$
and
 $$
  H_q^{(a-1)(t)}(\comega)\subset  \cases
 \dot
H_q^{t} (\overline\Omega )+e^+ d^{a-1}\overline H_q^{t-a+1}
                        (\Omega )\text{ when }t-a+1-\frac1q\notin \N,\\
  \dot
H_q^{t-\varepsilon } (\overline\Omega )+e^+ d^{a-1}\overline H_q^{t-a+1}
                        (\Omega )\text{ when }t-a+1-\frac1q\in \N,
                        \endcases
\tag3.8
$$

For the Dirichlet trace operator
 $\gamma _0^{a-1}u=\Gamma (a)\gamma _0(u/d^{a-1})$, defined on
 $H^{(a-1)(2a+s)}(\comega)$ when $a+s>\frac1q$, we have as a
 consequence of Theorem 3.2:

\proclaim{Corollary 3.4} 
Let $a+s>\frac1q$. Then $H^{(a-1)(2a+s)}(\comega)$ has a direct sum
decomposition
$$
H_q^{(a-1)(2a+s)}(\comega)=H_q^{a(2a+s)}(\comega)\; \dot+ \;
  d^{a-1}K_{(0)}B_q^{a+s+1-1/q}(\partial\Omega ).
\tag3.9
$$
Denote the elements $u=w+z$, then
 $w$ runs through $H_q^{a(2a+s)}(\comega)$ (equal to  $\dot
H_q^{2a+s}(\comega)$ when $a+s\in (\frac1q,1+\frac1q)$), and
$z=d^{a-1}K_{(0)}\varphi  $, with
$\varphi  $ running through $B^{a+s+1-1/q}_q(\partial\Omega )$. Moreover, $\varphi  =\gamma _0(u/d^{a-1})$.
\endproclaim

With the normalization (3.1), the formula for $z$ is
 $z=K^{a-1}_{(0)}\varphi ' $, $\varphi ' =\gamma ^{a-1}_0u$.
 
Also for higher values of $t$, there are direct sum decompositions where the first component is a supported
 space
 $\dot H^{t}(\comega)$ (or a transmission space); here the second component
is described by the composition of a vector of Poisson-type operators (derived from
$K_{(0)}$) with
 a vector of weighted traces.
They are covered by the following  general statement shown in \cite{G19}, Th.\
 3.10:

\proclaim{ Theorem 3.5} Let $\mu >-1$, $M\in{\Bbb N}$ and $t>\mu +M-\tfrac1{q'}$. Let
 $$
 \varrho ^\mu  _M = \{ \gamma ^\mu _0, \dots ,\gamma ^\mu 
_{M-1}\}.\tag3.10
$$
There is a vector of weighted Poisson-type operators $\widetilde{\Cal K}^\mu _M$ 
(cf.\ \cite{G19},
Theorem {\rm 3.9}), such that
 the elements $u\in H_q^{\mu (t)}(\comega)$ have unique decompositions 
$$
u=w+\widetilde{\Cal K}^\mu _M \cdot \varphi \in   H_q^{(\mu +M)(t)}(\comega)+e^+d^\mu \ol
H_q^{t-\mu }(\Omega )\cap H_q^{\mu (t)}(\comega),\tag3.11
$$
where 
$w $ runs through $ H_q^{(\mu +M)(t)}(\comega)$ (equal to $\dot
H_q^t(\comega)$ if $t-\mu \in(M-\tfrac1{q'},M+\tfrac1p )$), and
$\varphi   $ runs through $\prod _{0\le j<M}B_q^{t-\mu
-j-1/q}(\partial\Omega )$, when $u$ runs through $H_q^{\mu (t)}(\comega)$. Here
$\varrho  ^\mu _Mu$ equals  $\varphi   $.

The operator $\widetilde{\Cal K}^\mu _M$ satisfies that $\widetilde
{\Cal K}'_M=d^{-\mu
}\widetilde{\Cal K}^\mu _M$ is a vector
$$
\widetilde
{\Cal K}'_M=\{\widetilde
{\Cal K}'_{M,0},\dots \widetilde
{\Cal K}'_{M,M-1}\}
$$
of Poisson operators in the
Boutet de Monvel calculus; $\widetilde
{\Cal K}'_{M,j}$ being of order $-j$.
\endproclaim

The last statement
in the theorem refers to the description of the structure of $\widetilde{\Cal K}^\mu _M$ in
\cite{G19, (3.27)}, which contains additional information.

Versions of the
result already appeared in \cite{G15}, cf.\ Th.\ 5.4 there.

The Neumann trace operator  $\gamma _1^{a-1}u=\Gamma (a+1)\gamma
_1(u/d^{a-1})$ is included in the trace vector $\varrho
_2^{a-1}=\{\gamma _0^{a-1},\gamma _1^{a-1}\}$; for this we have as an
application of  Theorem 3.5 with $\mu =a-1$, $M=2$:

\proclaim{Corollary 3.6} Let $t>a+1-\frac1{q'}$. There is a 2-vector of
weighted Poisson-type operators $\widetilde{\Cal K}^{a-1} _2$ (where
 $d^{1-a}\widetilde{\Cal K}^{a-1} _2$ is a Poisson operator in
the Boutet de Monvel calculus),  such
that the elements $u\in H_q^{(a-1) (t)}(\comega)$ have unique decompositions 
$$
u=w+\widetilde{\Cal K}^{a-1} _2 \cdot \varphi \in
H_q^{(a+1)(t)}(\comega)+e^+d^{a-1} \ol
H_q^{t-a+1 }(\Omega )\cap H_q^{(a-1) (t)}(\comega);\tag3.12
$$
here 
$w $ runs through $ H_q^{(a+1)(t)}(\comega)$
and  $\varphi   $ runs through $\prod _{j=0,1}B_q^{t-a+1
-j-1/q}(\partial\Omega )$ when $u$ runs through $ H_q^{(a-1) (t)}(\comega)$. Moreover,
$\varrho  ^{a-1} _2u$ equals  $\varphi   $.

In particular, if $t \in(1+a-\tfrac1{q'},1+a+\tfrac1q )$, hence
$s=t-2a\in (1-a-\tfrac1{q'},1-a+\tfrac1{q})$, 
$$
H_q^{(a-1)(2a+s)}(\comega)=\dot H_q^{2a+s}(\comega)\; \dot+ \;
  \widetilde{\Cal K}^{a-1} _2{\prod }_{j=0,1}B_q^{a+s+1/q' -j}(\partial\Omega ).
\tag3.13
$$
With elements denoted $u=w+z$, $w$ runs through $\dot H_q^{2a+s}(\comega)$, and
$z$ equals $\widetilde{\Cal K}^{a-1} _2\cdot \varphi $, where $\varphi $ runs
through  $\prod _{j=0,1}B_q^{a+s+1/q' -j}(\partial\Omega )$; here   $\varphi =\varrho _2^{a-1}u$.

\endproclaim

\subsubhead 3.2 Results in  H\"older-Zygmund spaces 
\endsubsubhead

The decomposition can be set up in  $C^t_*$-spaces with a similar use
of $K_{(0)}$. Here we have from \cite{G19}, Th.\ 3.6:

\proclaim{Theorem 3.7}  With the Poisson operator  $K_{(0)}$
introduced in Theorem
{\rm 3.2}, and $K^\mu _{(0)}$ defined by {\rm (3.2)}, the following holds:

For $t>\mu $,  $K^{\mu }_{(0)}$ maps $C_*^{t-\mu }(\partial\Omega )\to
e^+d^{\mu }\ol C_*^{t-\mu }(\Omega )\cap C_*^{(\mu )(t)}(\comega)$.
Moreover, the elements of
 $C^{\mu (t)}_*(\comega)$ have unique decompositions
$$
u= K^\mu _{(0)}\varphi  +w, \text{ where }
\tag3.14$$
$\{\varphi ,w\}$ runs through $ C_*^{t-\mu }(\partial\Omega )\times
C_*^{(\mu +1)(t)}(\comega)$;   here $\varphi $ equals $\gamma ^\mu _0u$.

The space
$C_*^{(\mu +1)(t)}(\comega)$ equals $\dot C^t_*(\comega)$ if $t-\mu \in (0,1)$.
\endproclaim

Again there are corollaries in special cases:

Recall that the solution space for (2.4) when $f\in \ol C_*^s(\Omega )$
equals the $a$-transmission space $C_*^{a(2a+s)}(\comega)$. We have in
particular for
$a+s\in (0,1)$:

\proclaim{Corollary 3.8} Let 
 $s\in (-a,1- a)$.
The $a$-transmission 
space $C_*^{a(2a+s)}(\comega)$  
is a direct sum:
$$
  C_*^{a(2a+s)}(\comega)=\dot C_*^{2a+s}(\comega)\; \dot+ \;
  K^a_{(0)}C_*^{a+s}(\partial\Omega ).
\tag3.15
$$ 
With elements denoted $u=w+z$, $w$ runs through $\dot C_*^{2a+s}(\comega)$, and $z=K^a_{(0)}\psi
$, with
$\psi $ running through $C_*^{s+a}(\partial\Omega )$. Moreover, $\psi
=\gamma _0^au$.

\endproclaim

There are also results for higher values of $s$ and $t$. We have
from \cite{G19}, Th.\ 3.12:

\proclaim{Theorem 3.9} Let $M\in{\Bbb N}$.
The weighted Poisson-type operator $\widetilde{\Cal K}^\mu _M$
recalled in Theorem {\rm 3.5} maps 
$\prod _{0\le j<M}C_*^{t-\mu -j}(\partial\Omega )$
into  $  e^+d^{\mu }\ol C_*^{t-\mu }(\Omega )\cap  C_*^{\mu (t)}(\comega)$ when $t>\mu +M-1$.

 For  $t>\mu +M-1$,  the elements of
 $C^{\mu (t)}_*(\comega)$ have unique decompositions
$$
u=w+\widetilde{\Cal K}^\mu _{M}\cdot \varphi  \in  C_*^{(\mu 
+M)(t)}(\comega) +e^+d^\mu \ol
C_*^{t-\mu }(\Omega )\cap C_*^{\mu (t)}(\comega); 
\tag3.16$$
here $w$ 
runs through $
C_*^{(\mu +M)(t)}(\comega)$ (equal to  $\dot C^t_*(\comega)$ if $t-\mu \in (M-1,M)$), and  $\varphi =\varrho ^\mu _Mu$ runs
through $\prod _{j<M}C_*^{t-\mu  -j}(\partial\Omega )$.

Here $\widetilde{\Cal K}^\mu _M=d^\mu \widetilde {\Cal K}'_M$, where
$\widetilde {\Cal K}'_M$ is a vector of standard Poisson operators of
orders $\{0,-1,\dots,- M+1\}$, as recalled in Theorem {\rm 3.5}. 
\endproclaim

This allows a precise description of the solution spaces for (2.4) when $f$
belongs to higher-order H\"older-Zygmund spaces.
Further rules in H\"older-Zygmund spaces (like those above in
Bessel-potential spaces) are stated explicitly in \cite{G19}.

\subhead 4. Decomposition of spaces over  $C^{1+\tau} $-domains \endsubhead

For operators $P$ satisfying the $C^\tau $-version of Hypothesis 2.1,
it was shown in \cite{AG23} that the results on the description of the
Dirichlet domain as an $a$-transmission space can be generalized to
domains $\Omega $ with $C^{1+\tau }$-boundary; allowing the
smoothness parameter $t$ to lie in $[0,2a-\tau )$.  We shall
therefore investigate the generalization  of results in
Section 3
to $C^{1+\tau }$-domains.
The nonsmooth analysis was pursued further in \cite{G23}, where the
transmission spaces were analysed in more detail.

The definition of
$H_q^{\mu (t)}(\comega)$ 
when $\Omega $ is $C^{1+\tau }$ is formulated via local coordinate
systems, and is most convenient when $\tau \ge 1$, since the distance
function $d$  is then $C^{1+\tau }$ (see \cite{AG23, Sect.\ 2.1} or
\cite{G23, Sect.\ 2.1}). We take $\tau
\ge 1$ in the main applications.

Details of the definition are given in \cite{AG23, Sect.\ 4.3}, and also in
\cite{G23, Rem.\ 2.1} for the $H_q^t$-case.

\example{Remark 4.1} It is recalled in \cite{AG23, Sect.\ 4.2} that
a $C^1$-diffeomorphism $F\colon \rn\to\rn$ with $DF\in C^\tau
(\rn\times\rn)$  (some $\tau >0$) defines a mapping $F^*\colon
u\mapsto u\circ F$, that is bounded $H_q^t(\rn)\to H_q^t(\rn)$ when
$0\le t<1+\tau $. For the diffeomorphisms $F_\zeta \colon x\mapsto
(x',x_n-\zeta (x'))$ with $\zeta \in C^{1+\tau }(\R^{n-1})$ used in \cite{AG23}, \cite{G23} and the present
paper, the boundedness holds 
also when $t\in (-1-\tau ,0)$, by the duality argument explained in the
lines above (6.7) in \cite{AG23} (in the case $t=-a$). Negative $t$ (e.g.\ down
to $\mu -\frac1{q'}$) enter
in some formulations. 
\endexample

We begin by elaborating some properties set up in \cite{AG23, Sect.\ 4.2}.

It
is well-known that the multiplication by a function $\varphi \in
C^\sigma (\rn)$ is continuous in $H_q^t(\rn)$ and $C_*^t(\rn)$ when $|t|<\sigma $:
$$
\varphi H_q^t(\rn)\subset H_q^t(\rn),\quad \varphi C_*^t(\rn)\subset
C_*^t(\rn),\quad \text{ when }\varphi \in
C^\sigma (\rn), |t|<\sigma ;\tag4.1
$$
this is important for the localization procedure used to generalize the definitions of the spaces to domains $\Omega \subset
\rn$, and to transfer known properties for spaces
over $\rnp$ to spaces over
$\Omega $.

Multiplication by $\varphi $ is likewise continuous in $\mu $-transmission spaces under certain conditions on the parameters:

\proclaim{Proposition 4.2} Let $\mu >-1$ and $\sigma >0$, and let
$\varphi \in C^\sigma (\rn)$.

$1^\circ$ The localization property
$$
u\in H_q^{\mu (t)}(\crnp)
\implies \varphi u\in H_q^{\mu (t)}(\crnp)\tag 4.2
$$ holds under each of the
following conditions:
\roster
\item $
\mu \ge 0$, $\sigma >\mu $ and $\mu -\tfrac1{q'}<t<\sigma $,
\item $\mu \in (-1,0)$, $\sigma >\mu +1$ and $\mu -\tfrac1{q'}<t<\sigma -1$.
\endroster

$2^\circ$ The localization property
 $$
u\in C_*^{\mu (t)}(\crnp)
\implies \varphi u\in C_*^{\mu (t)}(\crnp)\tag 4.3
$$ holds under each of the
following conditions:
\roster
\item $t\in (\mu -1,\mu )$ and $\sigma \ge \max\{\mu ,1-\mu\}$,
\item $t \in [\mu ,\mu +1)$ and $t<\sigma$, where $
\sigma \ge  \mu_+$,
\item $t\ge \mu +1$, and $t<\sigma -1$ if $\mu \ge 1$,
$t<\sigma -2$ if $\mu \in [0,1)$, $t<\sigma -3$ if $\mu \in (-1,0)$.
\endroster

\endproclaim

\demo{Proof}
The case $1^\circ$ was proved in detail in \cite{AG23, Prop.\ 4.5}.
The case $2^\circ$ entered in \cite{AG23, Th.\ 4.6}, but was not
accounted for in detail there. We give a proof in Theorem A.2 in the
Appendix to the present paper. \qed
\enddemo

With this, we have as in  \cite{AG23, Sect.\ 4.2}:

\proclaim{Proposition 4.3} Let  $\mu >-1$ and $\tau >0$, and
let $\Omega $ be $C^{1+\tau }$.

$1^\circ$ The spaces $H_q^{\mu (t)}(\comega)$ are
defined for $\mu -\frac1{q'}<t<1+\tau $.

The
spaces satisfy {\rm 
(2.15)--(2.18)} under the additional conditions:
 If $\mu \ge 0$, assume $1+\tau >\mu $. If $\mu < 0$, assume $t<\tau $.

$2^\circ$ The spaces $C_*^{\mu (s)}(\comega)$ are defined for $\mu -1
<t<1+\tau $.

The spaces satisfy  {\rm (2.21)--(2.23)} under the additional
conditions in Proposition {\rm 4.2} $2^\circ$, with $\sigma $ replaced by $1+\tau $.
\endproclaim

\demo{Proof}  $1^\circ$. We briefly recall the {\it localization procedure} from \cite{AG23, Def.\
4.3}: When $\zeta \in C^{1+\tau }(\R^{n-1},\R)$, defining the curved
half-space $\rn_\zeta =\{x\in\rn \mid x_n>\zeta (x')\}$, $
H_q^{\mu (t)}(\ol {\R}^n_\zeta )$ is defined as $F^*_\zeta ( H_q^{\mu
(t)}(\crnp))$ using the diffeomorphism $F_\zeta \colon x \mapsto
(x',x_n-\zeta (x'))$; here 
$$
 u\in 
H_q^{\mu (t)}(\ol {\R}^n_\zeta  )\iff u\circ F_\zeta ^{-1}\in  H_q^{\mu
(t)}(\crnp).\tag4.4
$$
For a bounded $C^{1+\tau }$-domain $\Omega $, each point $x_0\in \partial\Omega $
has a bounded open neighborhood $U\subset \rn$ and a  $\zeta \in
C^{1+\tau }(\R^{n-1},\R)$ such that (after a suitable rotation)
$\Omega \cap U=\rn_\zeta \cap U$. The space $H_q^{\mu (t)}(\comega)$
consists of the functions $u\in H^t_{q,loc}(\Omega )$ such that for
each $x_0\in \partial\Omega $, with a $\varphi \in C_0^\infty (U)$
that is 1 on a neighborhood of $x_0$, $(\varphi u)\circ F_\zeta
^{-1}\in H_q^{\mu
(t)}(\crnp)$. (This is the same procedure as is used to define other
function spaces over $\comega$, e.g.\ $\ol H_q^t(\Omega )$.)

The validity of the properties {\rm 
(2.15)--(2.16)} is accounted for in \cite{AG23, Th.\ 4.6} (with preparations in
Prop.\ 4.4 and 4.5 there).
Note that the pull-back of
$\varphi $ to $\crnp$ is only $C^{\tau +1}$; the rule in Proposition
4.2 $1^\circ$ is
needed here.
The property (2.18) follows likewise from the result in smooth cases by localization.

$2^\circ$. This follows in the same way from a $C_*$-version of  \cite{AG23,
Def.\ 4.3 -- Th.\ 4.6}, relying on Proposition 4.2 $2^\circ$. \qed
\enddemo

We note in particular that  \cite{AG23,Prop.\
4.5} (with $\sigma =1+\tau $) requires when $\mu \ge 0$ that $1+\tau
>\mu $ and that $t<1+\tau $. For $\mu <0$ it requires
that $t<\tau $.

The fact that $t<\tau $ is only needed for $\mu \in (-1,0)$ was not
 pointed out explicitly in \cite{AG23, Th.\ 4.6}, but we will take it
 into account in the present paper.

 Further comments on the
localization are given in \cite{G23, Rem.\ 2.1}.

In a description of the finer structure, we shall use Prop.\ 2.2 and Thm.\ 2.3
from \cite{G23}. Prop.\ 2.2 shows the existence of a continuous
right inverse of $\gamma _0^\mu $, but we need to make this description
more explicit, so we now give an elaborated
version of Prop.\ 2.2. Here we rely on the extension of the Boutet de Monvel
calculus to nonsmooth symbols defined in \cite{A05,AG23}.

\proclaim{Theorem 4.4} Let  $\mu >-1$, $\tau >0$, and $\mu +\frac1{q}< t<1+\tau $ with
  $t-\mu<1+\tau$. 
 Let $\rn_\zeta = \{x\in\R^n\mid x_n>\zeta (x')\}$ be defined by a
  function $\zeta \in C^{1+\tau}(\R^{n-1})$, and let $d_1$ be  a
  bounded positive $C^{1+\tau }$-function satisfying $d_1(x)=x_n-\zeta (x')$ near
  $\partial\rn_\zeta $.

The Poisson operator $K_0=\operatorname{OPK}((\ang{\xi '}+i\xi
_n)^{-1})$ on $\rnp$ is transformed by the diffeomorphism  $F_\zeta \colon
(x',x_n)\mapsto (x',x_n-\zeta (x'))$ to an operator $K_{(0),\zeta }$
(a $C^\tau $-Poisson operator as in \cite{A05} plus a remainder),
which is a right inverse for $\gamma _0$ on $\rn_\zeta $. Moreover:

The mapping $\gamma
  _0^{\mu,1} \colon u\mapsto \Gamma (\mu +1)(u/d_1^\mu )|_{\partial
           \rn_\zeta }$
  is continuous and surjective:
$$
\gamma _0^{\mu ,1}\colon H_q^{{\mu} (t)}(\ol{\R}^n_\zeta )\to B_{q}^{t-\mu
  -\frac1q}(\partial \rn_\zeta),\tag4.5
$$
having as a continuous right inverse the operator $K_{(0),\zeta }^{\mu
,1}=\tfrac1{\Gamma (\mu +1)}d_1^\mu K_{(0),\zeta }$, mapping for all $t$:
$$
K_{(0),\zeta }^{\mu ,1}\colon B_{q}^{t-\mu
  -\frac1q}(\partial \rn_\zeta) \to H_q^{{\mu} (t)}(\ol{\R}^n_\zeta
  )\cap d_1^\mu \ol H^{t-\mu }(\rn_\zeta ).\tag4.6
$$
 Moreover, the space $H_q^{{(\mu+1)}
(t)}(\ol{\R}^n_\zeta  )$ is a closed subspace of $H_q^{{\mu} (t)}(\ol{\R}^n_\zeta  )$,
equal to the kernel of the mapping {\rm (4.5)}.

When $\tau \ge 1$, we can replace $d_1$ in {\rm (4.5)--(4.6)} by a
distance function $d$ (a $C^{1+\tau }$-function
equal to
$\operatorname{dist}(x,\partial\Omega )$ near $\partial\Omega $),
defining the weighted trace $\gamma _0^\mu \colon u\mapsto \Gamma (\mu
+1)\gamma _0(u/d^\mu )$ and having the right inverse $K_{(0),\zeta }^\mu =\tfrac1{\Gamma (\mu +1)}d^\mu K_{(0),\zeta }$,
with the same properties as listed  above for $K_{(0),\zeta }^{\mu
,1}$.

\endproclaim

\demo{Proof} The case where 
$\zeta (x')\equiv 0$
(the flat case), was recalled above in Lemma 3.1 and the subsequent statements
for $\Omega =\rnp$.
They carry over to
the case of general $\zeta (x')$ in view of the  mapping properties of the
diffeomorphism $F_\zeta \colon (x',x_n)\mapsto (x',x_n-\zeta (x'))$ and the
definitions 
listed in \cite{AG23, Sect.\ 4.2}. Note that the space $H_q^{{(\mu+1)}
(t)}(\ol{\R}^n_\zeta  )$ is well-defined with the parameter $\mu +1$, since
the hypotheses assure that  $t>(\mu +1)-1/{q'}$.

The operator $K_{(0),\zeta }$ arising from application of the
 diffeomorphism
 $F_\zeta $ can be seen to have the structure of a
Poisson operator with $C^\tau $ $x$-dependent symbol (as in
\cite{A05}) plus a remainder, shown as in Theorem 1.2.$2^\circ$
in \cite{AG23}. We refrain from a further analysis since we shall just
use the established mapping properties.

In the flat case,
$K_0^\mu
=\frac1{\Gamma (\mu +1)}x_n^\mu e^+K_0$ acts as described in Lemma 3.1. This is transformed to
$K_{(0),\zeta }^{\mu ,1}$ by the diffeomorphism $F_\zeta \colon
(x',x_n)\mapsto (x',x_n-\zeta (x'))$, which also transforms $x_n^\mu $
to $d_1^\mu $ (near the boundary); and (4.5)--(4.6) hold.

Finally, when $\tau \ge 1$, the
true distance function $d$ is $C^{1+\tau }$ with $g=d_1/d$ and $d/d_1$
being $C^\tau $, so an insertion of $d_1=gd$ leads to the definition
of $K^\mu _{(0),\zeta }$, likewise satifying the properties listed for $K^{\mu ,1}_{(0),\zeta }$.
\qed  
\enddemo

Note that $K^\mu _{(0),\zeta }$ is constructed independently of the parameter
$t$ in the function spaces.

A completely analogous result holds for operators over $C^t_*$-spaces. We
leave out the formulation for the most general $\tau >0$ (which can be invoked
whenever needed), and just formulate the case $\tau \ge 1$:

\proclaim{Theorem 4.5} Let $\mu >-1$, $\tau \ge 1$, and $\mu < t<1+\tau $ with
  $t-\mu<1+\tau$.

Let $\rn_\zeta = \{x\in\R^n\mid x_n>\zeta (x')\}$ be defined by a
  function $\zeta \in C^{1+\tau}(\R^{n-1})$, and let $d$ be  a
  distance function, as in Theorem {\rm 4.4}.
  
The mapping $\gamma
  _0^{\mu} \colon u\mapsto \Gamma (\mu +1)(u/d^\mu )|_{\partial
           \rn_\zeta }$
  is continuous and surjective:
$$
\gamma _0^{\mu }\colon C_*^{\mu (t)}(\ol{\R}^n_\zeta )\to C_*^{t-\mu
 }(\partial \rn_\zeta),\tag4.7
$$
having the continuous right inverse $K_{(0),\zeta }^\mu $ introduced
in 
Theorem {\rm 4.4},
mapping for all $t$:
$$
K_{(0),\zeta }^{\mu }\colon C_*^{t-\mu }
 (\partial \rn_\zeta) \to C_*^{{\mu} (t)}(\ol{\R}^n_\zeta
  )\cap d^\mu \ol C_*^{t-\mu }(\rn_\zeta ).\tag4.8
$$
 Moreover, the space $C_*^{{(\mu+1)}
(t)}(\ol{\R}^n_\zeta  )$ is a closed subspace of $C_*^{{\mu} (t)}(\ol{\R}^n_\zeta  )$,
equal to the kernel of the mapping {\rm (4.7)}.
\endproclaim

\demo{Proof} The proof of Theorem 4.4 extends to this setting
directly, since the mapping properties of $\gamma _0$, $K_0$ and $K^\mu _0$ hold in
the $C_*$-scale of spaces. \qed
\enddemo

The construction is carried over to $\Omega $ by localization as in
Proposition 4.3 (also explained in  \cite{G23, Th.\ 2.3}). It is here that the
multiplication properties in Proposition 4.2 are needed. We now recall the theorem and
provide it with more details.

\proclaim{Theorem 4.6} Let $\mu >-1$, $\tau \ge 1$, $1+\tau >\mu $ and $\mu
 +\frac1{q}< t<\tau +1$,
 and  let $\Omega\subset\R^n$ be a bounded
  $C^{1+\tau}$-domain. If $\mu \in (-1,0)$, assume moreover $\tau >\mu
  $ and  $t-\mu <\tau $.

  The mapping $\gamma _0^\mu \colon u\mapsto \Gamma (\mu +1)(u/d^\mu
)|_{\partial \Omega }$ is continuous and surjective:
$$
\gamma _0^\mu \colon H_q^{{\mu} (t)}(\comega )\to B_q^{t-\mu
  -\frac1q}(\partial \Omega ),\tag4.9
$$
having a  continuous right inverse $K^{\mu }_{(0)}$, mapping for all $t$
$$
K_{(0)}^{\mu }\colon B_{q}^{t-\mu
  -\frac1q}(\partial \Omega ) \to H_q^{{\mu} (t)}(\comega
  )\cap d^\mu \ol H_q^{t-\mu }(\Omega  ).\tag4.10
$$
 Moreover, the space $H_q^{{(\mu +1)}
(t)}(\comega )$ is a closed subspace of $H_q^{{\mu} (t)}(\comega )$,
and equals the kernel of the mapping {\rm (4.9)};
$$
 \{u\in H_q^{\mu (t)}(\comega)\mid\gamma
_0^{\mu }u=0\}=H_q^{(\mu +1)(t)}(\comega). \tag4.11
$$

\endproclaim

\demo{Proof} The continuity of the mapping $\gamma _0^\mu $ is
established in \cite{AG23 Th.\ 4.5} by use of a cover
$\bigcup_{i=0,1,\dots, I}U'_1$ (with $\overline U'_0\subset \Omega $) and an associated partition of unity
$\{\varrho _i\}_{i=0,\dots,I}$ , such that the $U'_i$ with $i\ge 1$
cover $\partial\Omega $ and for each such $U'_i$ there is a function
$\zeta _i\in C^{1+\tau }(\R^{n-1},\R)$  such that, after a
rotation and translation depending on $i$, $\Omega \cap U'_i =
\rn_{\zeta _i}\cap U'_i$. In each such neighborhood $U'_i$, the facts
known for $\rn_{\zeta _i}$ can be applied to $\varrho _iu$, this is
collected to a statement on $u$ by summation. (The conditions on
parameters assure that multiplication by cutoff functions is continuous.)

More precisely, write $u=\sum_{0\le i\le I}\varrho _iu$, whereby $$
\varphi =\gamma _0^\mu u =\sum_{1\le i\le I}\gamma _0^\mu (\varrho
_i u)=\sum_{1\le i\le I}\varphi _i,\; \varphi _i=\gamma _0^\mu (\varrho _iu).
$$
Let $\psi _i\in C_0^\infty (U_i)$, equal to 1 on a neighborhood of
$\supp \varrho _i$. In each $U_i$ we apply Theorem 4.4 (in the
translated and rotated situation; this is tacitly understood in the notation) to
construct a Poisson-type operator $K^\mu _{(0),\zeta _j}$ as a right
inverse of $\gamma _0^\mu $ with the properties in Theorem
4.4. Now let 
$$
K^\mu _{(0)}=\sum_{1\le i\le I}\psi _iK^\mu _{(0),\zeta _j}\varrho _j.
$$
It has the asserted mapping properties, and is a right inverse of
$\gamma _0^\mu $, since
$$
\aligned
\gamma _0^\mu K^\mu _{(0)}\varphi& =\sum_{1\le i\le I}\gamma _0^\mu
(\psi _iK^\mu _{(0),\zeta _j}\varrho _j\varphi )\\
&=\sum_{1\le i\le
I}\gamma _0 (\psi _i)\gamma _0^\mu (K^\mu _{(0),\zeta _j}\varrho
_j\varphi )=\sum_{1\le i\le I}\gamma _0 (\psi _i)
\varrho _j\varphi=\varphi ;\endaligned
$$
using that $\psi _i$ is 1 on $\supp \varrho _i$. Since we have for
each piece $\varrho _iu$ that 
$$
\varrho _iu-\psi _iK^\mu _{(0),\zeta _i}\gamma _0^\mu (\varrho _iu)=
\psi _i\varrho _iu-\psi _iK^\mu _{(0),\zeta _i}\gamma _0^\mu (\varrho _iu)
\in H_q^{(\mu +1)(s)}(\rn_{\zeta _i}),
$$
it follows by summation over $i$ (using also that $\varrho _0u\in \dot
H_q^s(\comega)$), that 
$$
u-K^\mu _{(0)}\gamma _0^\mu u\in  H_q^{(\mu +1)(s)}(\comega),
$$
verifying (4.11).
\qed
\enddemo

The  proof can be adapted directly to the $C^t_*$-setting, and gives:

\proclaim{Theorem 4.7} Let $\mu >-1$, $\tau \ge 1$, $t>\mu $,
and  let $\Omega\subset\R^n$ be a bounded
  $C^{1+\tau}$-domain. Assume that one of the conditions in Proposition
  {\rm 4.2 $2^\circ$} holds with $\sigma =1+\tau $.

  The mapping $\gamma _0^\mu \colon u\mapsto \Gamma (\mu +1)(u/d^\mu
)|_{\partial \Omega }$ is continuous and surjective:
$$
\gamma _0^\mu \colon C_*^{{\mu} (t)}(\comega )\to C_*^{t-\mu
}(\partial \Omega ),\tag4.12
$$
and $K^{\mu }_{(0)}$ defined in Theorem {\rm 4.6} acts as a right inverse, mapping for all $t$
$$
K_{(0)}^{\mu }\colon C_*^{t-\mu
}(\partial \Omega ) \to C_*^{{\mu} (t)}(\comega
  )\cap d^\mu \ol C_*^{t-\mu }(\Omega  ).\tag4.13
$$
 Moreover, the space $C_*^{{(\mu +1)}
(t)}(\comega )$ is a closed subspace of $C_*^{{\mu} (s)}(\comega )$,
and equals the kernel of the mapping {\rm (4.12)};
$$
 \{u\in C_*^{\mu (t)}(\comega)\mid\gamma
_0^{\mu }u=0\}=C_*^{(\mu +1)(t)}(\comega). \tag4.14
$$

\endproclaim

The theorems show in particular for $\mu =a$:

\proclaim{Corollary 4.8} Let $a\in (0,1)$ and $\tau \ge 1$, and  let $\Omega\subset\R^n$ be a bounded
  $C^{1+\tau}$-domain.

$1^\circ$  When $a +\frac1{q}< t<\tau +1 $,
 the mapping $\gamma _0^a \colon u\mapsto \Gamma (a +1)(u/d^a
)|_{\partial \Omega }$ is continuous and surjective:
$$
\gamma _0^a \colon H_q^{{a} (t)}(\comega )\to B_q^{t-a
  -1/q}(\partial \Omega ),\tag4.15
$$
having a  continuous right inverse $K^{a }_{(0)}$. Moreover, the space $H_q^{{(a +1)}
(t)}(\comega )$ is a closed subspace of $H_q^{{a} (t)}(\comega )$,
and equals the kernel of the mapping {\rm (4.15)};
$$
 \{u\in H_q^{a (t)}(\comega)\mid\gamma
_0^{a }u=0\}=H_q^{(a+1)(t)}(\comega). \tag4.16
$$

$2^\circ$
Let either $a < t< a+1$, or  $a+1\le t<\tau -1$. Then  the mapping $\gamma _0^a \colon u\mapsto \Gamma (a +1)(u/d^a
)|_{\partial \Omega }$ is continuous and surjective:
$$
\gamma _0^a \colon C_*^{{a} (t)}(\comega )\to C_*^{t-a
}(\partial \Omega ),\tag4.17
$$
having a  continuous right inverse $K^{a }_{(0)}$ (consistent with the
operator in $1^\circ$). Moreover, the space $C_*^{{(a +1)}
(t)}(\comega )$ is a closed subspace of $C_*^{{a} (t)}(\comega )$,
and equals the kernel of the mapping {\rm (4.17)};
$$
 \{u\in C_*^{a (t)}(\comega)\mid\gamma
_0^{a }u=0\}=C_*^{(a+1)(t)}(\comega). \tag4.18
$$

\endproclaim

For $t=2a+s$ with $s$ suitably close to 0, this gives the following
generalizations of Corollary 3.3 and 3.8:

\proclaim{Corollary 4.9} Let $a\in (0,1)$ and $\tau \ge 1$, and  let $\Omega\subset\R^n$ be a bounded
  $C^{1+\tau}$-domain.

$1^\circ$ For $s\in (-a+\frac1q,-a+1+\frac1q)$ with $s<\tau +1-2a$,
$H_q^{a(2a+s)}(\comega)$
is a direct sum:
$$
H_q^{a(2a+s)}(\comega)=\dot H_q^{2a+s}(\comega)\; \dot+ \;
  K^a_{(0)}B_q^{a+s-1/q}(\partial\Omega ).
\tag4.19
$$
With elements denoted $u=w+z$, $w$ runs through $\dot H_q^{2a+s}(\comega)$, and
$z$ equals $K^a_{(0)}\psi $, with
$\psi $ running through $B^{a+s-1/q}_q(\partial\Omega )$. Moreover,
$\psi =\gamma ^a _0u$.

$2^\circ$ For $s\in (-a,1-a)$, $C_*^{a(2a+s)}(\comega)$
is a direct sum:
$$
C_*^{a(2a+s)}(\comega)=\dot C_*^{2a+s}(\comega)\; \dot+ \;
  K^a_{(0)}C_*^{a+s}(\partial\Omega ).
\tag4.20
$$
With elements denoted $u=w+z$, $w$ runs through $\dot C_*^{2a+s}(\comega)$, and
$z$ equals $K^a_{(0)}\psi $, with
$\psi $ running through $C_*^{a+s}(\partial\Omega )$. Moreover,
$\psi =\gamma ^a _0u$. 
\endproclaim

\demo{Proof} $1^\circ$. The space $H_q^{(a+1)(t)}(\comega)$ equals
$\dot H_q^t(\comega)$, when $t\in (a+\frac1q,a+1+\frac1q)$ (besides
satisfying $t<\tau +1$). Then the
decomposition follows from Corollary 4.8 $1^\circ$ since $K^a_{(0)}$ is a right inverse of
$\gamma _0^a$ with the asserted mapping properties. The result is
stated in terms of $s=t-2a$, and the conditions on $t$ carry over to
the conditions $s\in (-a+\frac1q,-a+1+\frac1q)$, $s<\tau +1-2a$.

$2^\circ$. The space $C_*^{(a+1)(t)}(\comega)$ equals
$\dot C_*^t(\comega)$, when $t\in (a,a+1)$. Here the
decomposition follows from Corollary 4.8 $2^\circ$ since $K^a_{(0)}$ is a right inverse of
$\gamma _0^a$ with the asserted mapping properties. In terms of
$s=t-2a$,  the conditions on $t$ translates to $s\in (-a,1-a)$.
\qed

\enddemo

Note that the case $s=0$ enters in $2^\circ$, and in $1^\circ$ when
$q>1/a$; in both cases $\tau =1$ is allowed.

There are also generalizations of the corollaries in case $\mu =a-1$
and in cases considering vectors of weighted boundary values, as in
Corollary 3.4, Theorem 3.5, Corollary.3.6, Theorem 3.9; they can be
invoked in a similar way whenever necessary.

\subhead 5. Embeddings with powers of the boundary distance \endsubhead

For further developments of the above results we include a study of
how the boundary distance enters in embeddings of spaces $\dot C_*^{t
}(\comega )$ and $\dot H_q^{t}(\comega )$.
First there is a result for H\"older spaces:

\proclaim{Proposition 5.1} Let $\Omega $ be a bounded $C^{1+\tau }$-domain, $\tau
>0$.
Let $\alpha $ and $\beta >0$ with $\alpha +\beta <1+\tau $. For $\alpha ,\beta ,\alpha +\beta \notin{\Bbb N}$,
there holds
$$
\dot C^{\alpha +\beta }(\comega)\subset d(x)^\alpha  \dot C^{\beta }(\comega ).\tag5.1
$$

The inclusion also holds when $\Omega $ is a curved halfspace $\R^n_\zeta
=\{x=(x',x_n)\mid x_n>\zeta (x')\}$, defined from a function
$\zeta \in C^{1+\tau }(\R^{n-1})$, and provided with the distance
function $d_1(x)=x_n-\zeta (x')$ for $x_n\le K+1$ (some  $K\ge \sup
|\zeta |$), extended smoothly to the rest of $\rnp$ as a positive
constant for large $x_n$.
\endproclaim

The result is apparently well-known, at least for bounded sets. For a
detailed proof see \cite{G23, Lemma A.5, Lemma A.6}.

For the scale of Bessel-potential spaces $\dot H^t_p(\comega )$, the
following result is known:

\proclaim{Proposition 5.2} Let  $p\in (1,\infty )$ and $t>0$.

$1^\circ$ 
For $f\in \dot
H^t_p(\crnp)$ with $\supp f\subset \ol B=\{|x|\le 1\}$
there is a constant $C$ 
such that
$$
\|x_n^{-t}f\|_{L_p(\rnp)}\le C\|f\|_{\dot H^t_p(\crnp)}.\tag5.2
$$

$2^\circ$ Let $\tau \ge 1$, $0<t<1+\tau $, and let $\Omega \subset \rn$ be
bounded and 
$C^{1+\tau }$. There is a constant $C$ such that for all $f\in \dot H^t_p(\comega)$,
$$
\|d^{-t}f\|_{L_p(\Omega )}\le C\|f\|_{\dot H^t_p(\comega)}.\tag5.3
$$

\endproclaim

\demo{Proof}
$1^\circ$.
This is shown in Triebel
 \cite{T01, Prop.\ 5.7} and its proof, applied to the scale
 $H^t_{p}(\rn)=F^t_{p2}(\rn)$ (in fact valid also for many other scales of function spaces).

$2^\circ$. Triebel concludes this for $C^\infty $-domains. To allow
$C^{1+\tau }$-domains, we first note that the estimate (5.2) carries over to
the case of curved half-spaces $\R^n_\zeta $ with $\zeta \in C^{1+\tau
}(\R^{n-1})$ as described in Theorem 4.4, since the diffeomorphism
$F_\zeta $ preserves $H^t_p(\rn)$ (for $|t|<1+\tau $), with $x_n$
carried over to $d_1(x)=x_n-\zeta (x')$. Since $\tau \ge 1$, $d_1$ can
be replaced by $d$ over bounded sets, as in Theorem 4.4. Then this is used
to show the
property for bounded $C^{1+\tau }$-domains, by localization as in the
proof of Theorem 4.6.
\qed

\enddemo

One can of course replace the ball $B=\{|x|<1\}$ by any other ball
$B_R=\{|x|<R\}$, $R>0$.

We shall now show that embeddings similar to those in Proposition 5.1
can be obtained for higher-order Bessel-potential spaces. To our
knowledge, this is a new result.
Here we first treat the
case with integer powers of $x_n$. (The subscript $p$ is used instead
of $q$ in this section, to facilitate the references to \cite{T01}.)

\proclaim{Proposition 5.3} Let $p\in (1,\infty )$. For $t\in [0,\infty )$, $s\in \N_0$, there
holds  when $\supp u\subset \ol B $,
$$
u\in \dot H_p^{s+t}(\crnp)\implies  u\in x_n^s\dot H_p^{t}(\crnp);\tag5.4
$$
in other words 
$$
\|x_n^{-s}u\|_{\dot H_p^{t}(\crnp)}\le C \|u\|_{\dot H^{s+t}_p(\crnp)}.\tag5.5
$$

\endproclaim

\demo{Proof} In these statements, $x_n^{-0}$ is read as 1.
Recall that for $k\in \N_0$, $H_p^k(\rn)$ is the classical
Sobolev-Slobodetskii space consisting of functions $u\in L_p(\rn)$
with all partial derivatives up to order $k$ in $L_p(\rn)$; complex
interpolation gives the noninteger cases. Note also
that to show
that $u\in H_p^r(\rn)$ for some $r\in [1,\infty )$, it suffices to show that
$\partial_1u,\dots, \partial_nu$ and $u$ are in
$H_p^{r-1}(\rn)$. Recall that $\supp u\subset \ol B$ throughout the proof.

First assume that $t$ is integer.
The calculations will
be made for increasing values of $k=s+t$. The cases where
$s=0$ are trivial (they are just included for ease of formulation). The
cases where $t=0$ are covered directly by Proposition 5.2.

{\it The case $k=1$.} The validity of (5.4) is secured by Proposition 5.2:
$$
u\in \dot H_p^1(\crnp)\implies u\in x_nL_p(\rnp).
$$

Proceed by induction: Assume that the embeddings (5.4) (with (5.5)) have been proved
for $s+t$ up to a value $k$, and consider the problem for $s+t=k+1$.

{\it The case $k+1$}.  Let $u\in \dot H_p^{k+1}(\crnp)$. We consider
successively the subcases where $t=0,1,\dots,k$. First, the case
$t=0$, $s=k+1$, follows obviously from
Proposition 5.2. In the next case $t=1$,  $s=k$,  we must show that
$$
x_n^{-k}u\in \dot H_p^1(\crnp).\tag5.6
$$
Here
$$
\partial_n(x_n^{-k}u)=-k x_n^{-k-1}u+x_n^{-k}\partial_nu. 
\tag5.7
$$
The first term is in $L_p(\rnp)$ by Proposition 5.2. For the second term,
$\partial_nu\in \dot H_p^k(\crnp)$ implies by Proposition 5.2 that the term
is in $L_p(\rnp)$. Since 
$\partial_j(x_n^{-k}u)=x_n^{-k}\partial_ju$  for $j<n$, Proposition
5.2 likewise gives that  $\partial_j(x_n^{-k}u)\in L_p(\rnp)$. Also
$x_n^{-k}u\in L_p(\rnp)$, so it follows that $x_n^{-k}u\in \dot H_p^1(\crnp)$, and
(5.6) has been obtained.

For higher $t$, we shall show how the validity in the case
$t=m+1,s=k-m$ follows from the validity when $t=m,s=k-m+1$ for each
$m$, where $m$ runs through 1 up to $k$. Here we use the formulas
$$
\aligned
\partial_n(x_n^{m-k}u)&=(m-k)x_n^{m-k-1}u+x_n^{m-k}\partial_nu,\\
\partial_j(x_n^{m-k}u)&=x_n^{m-k}\partial_ju,\text{ for }j<n.
\endaligned
\tag5.8
$$
Assume that we already know that $x_n^{m-1-k}u\in \dot H_p^m(\rnp)$.
Since $\partial_nu\in \dot H_p^k(\rnp)$, where the induction hypothesis
assures that $x_n^{m-k}\partial_nu\in \dot H_p^m(\rnp)$, the first
line in (5.8) shows that $\partial_n(x_n^{m-k}u)\in \dot
H_p^m(\rnp)$. The second line shows in view of the induction
hypothesis that  $\partial_j(x_n^{m-k}u)\in \dot
H_p^m(\rnp)$ for $j<n$. It follows that $x_n^{m-k}u\in \dot
H_p^{m+1}(\rnp)$, as was to be shown.
This ends the proof that (5.4) (and (5.5)) holds for $s+t=k+1$,
completing the induction step frem $k$ to $k+1$. It follows that (5.4)
holds for all $s,t\in\N_0$.

We now have the validity of (5.4) and (5.5) when $s$ and $t$ are
integers.
To include noninteger values of $t$,  fix $s\in \N$. We know that the
multiplication mappings $x_n^{-s}$:
$$
u\in \{w\in\dot H_p^{s+t}(\crnp)\mid \supp w\subset \ol B \} \mapsto
x_n^{-s}u\in \{w\in\dot H_p^{t}(\crnp)\mid \supp w\subset \ol B \} \tag5.9
$$
are continuous for $t\in \N_0$. By complex interpolation, the continuity extends
to noninteger values of $t$, allowing $t\in [0,\infty )$. \qed
\enddemo

Next, we show the full result including noninteger powers of $x_n$:

\proclaim{Theorem 5.4} Let $p\in (1,\infty )$.

$1^\circ$ For functions $u\in e^+L_p(\rnp)$
 with
 $\supp u\subset \ol B $,
{\rm (5.4)} and  {\rm (5.5)} hold for all $s,t\in \crp$.

$2^\circ$  Let $\Omega $ be a bounded $C^{1+\tau }$-domain, $\tau
\ge 1$, and let $s,t\in \crp$ with $s+t<\tau $. Then there holds
$$
u\in \dot H_p^{s+t}(\comega)\implies  u\in d^s\dot H_p^{t}(\comega);\tag5.10
$$
in other words 
$$
\|d^{-s}u\|_{\dot H_p^{t}(\comega)}\le C \|u\|_{\dot H^{s+t}_p(\comega)}.\tag5.11
$$
For $\tau \in (0,1)$ there holds a local version of the embeddings,
inferred from {\rm (5.13)} below.
\endproclaim

\demo{Proof} $1^\circ$. Fix $s\in \rp$. For integer values
$t=0,1,2,\dots$, we  show the validity of
the statements by induction. Recall that $\supp u\subset \ol B$ is
assumed throughout.

For $t=0$, the statements hold by Proposition 5.2.

Next, let $t=1$. Let $u\in
 \dot H_p^{s+1}(\crnp)$. To check whether $x_n^{-s}u\in 
\dot H_p^{1}(\crnp)$, write:
$$
\partial_n(x_n^{-s}u)=x_n^{-s}\partial_nu -sx_n^{-1-s}u = x_n^{-s}(\partial_nu-sx_n^{-1}u)
.\tag5.12
$$
Here  $\partial_nu\in \dot H_p^{s}(\crnp)$.
Moreover, since $u\in \dot H_p^{s+1}(\crnp)$, $x_n^{-1}u\in \dot
H_p^s(\crnp)$ by Proposition 5.3, and is likewise
supported in $\ol B $. So the parenthesis in (5.12) is in $\dot
H_p^{s}(\crnp)$, and then by Proposition 5.2, the full expression is
in $L_p(\rnp)$.
It is easy to check that also the tangential
 derivatives  $\partial_j(x_n^{-s}u)$  ($j<n$) and $x_n^{-s}u$
 itself are in  $L_p(\rnp)$, so
 we can conclude that $x_n^{-s}u\in \dot H_p^1(\crnp)$. This shows the
 statement for $t=1$.

The induction step: Assume that (5.4), (5.5) have been shown for
integers $t$ up to $k$, and let  $u\in
 \dot H_p^{s+k+1}(\crnp)$. To check whether $x_n^{-s}u\in \dot
H_p^{k+1}(\crnp)$, we use again the calculation (5.12).
Clearly,  $\partial_nu\in \dot H_p^{s+k}(\crnp)$, and $x_n^{-1}u$ is
there by Proposition 5.3. Then by the case $t=k$,
$x_n^{-s}(\partial_nu-sx_n^{-1}u)\in \dot H_p^{k}(\crnp)$.
With similar
estimates for the tangential derivatives $\partial_j(x_n^{-s}u)$, we conclude
that $x_n^{-s}u\in \dot H^{k+1}_p(\crnp)$, as was to be shown.

We have then obtained the validity of (5.4), (5.5) for all integer
$t$. This means that the mappings (5.9) are continuous for all $t\in
\N_0$, and it follows by complex interpolation that they are
continuous for all $t\in \crp$. This ends the proof of $1^\circ$.

$2^\circ$. To deduce the result for curved sets $\Omega $ we use the
tools recalled in Section 4, see also \cite{AG23, Sect.\ 4.2}. By a diffeomorphism $F_\zeta \colon x\mapsto
(x',x_n-\zeta (x'))$ with $\zeta \in C^{1+\tau }(\R^{n-1})$,
$\dot H_p^r(\crnp)$ is carried over to  $\dot H_p^r(\ol{\R}^n_\zeta)$
when $0\le r<1+\tau $, $x_n$ being carried over to the distance function
$d_1(x)=x_n-\zeta (x')$ near $\partial \ol{\R}^n_\zeta$. Then the
embedding (5.4) implies for $u$ supported in the image $\tilde B$ of  $B $:
$$
u\in \dot H_p^{s+t}(\ol{\R}^n_\zeta)\implies  u\in d_1^s\dot
H_p^{t}(\ol{\R}^n_\zeta).\tag5.13
$$
This gives a local information pertaining to the embedding (5.10). If
$\tau \ge 1$, we can replace $d_1$ by the true distance $d$ in open
neighborhoods $U$ of points of 
$\partial\Omega $ where  $r_U\dot
H_p^{s+t}(\ol{\R}^n_\zeta)$ and  $r_Ud_1^sr_{\tilde B}\dot
H_p^{t}(\ol{\R}^n_\zeta)$ are considered, so that in the description of the spaces in
local coordinates, the information (5.13) applied to such neighborhoods leads to (5.10).  \qed
 
\enddemo

To the best of our knowledge, the general result of Theorem 5.4 is new.

As a simple
consequence, we can also address multiplication by positive powers of $d$:

\proclaim{Theorem 5.5} Let $p\in (1,\infty )$. Let $s,t\in \crp$.

$1^\circ$ There holds for $u$ supported in $\ol B$:
$$
\|x_n^{s}u\|_{\dot H_p^{t}(\crnp)}\le C \|u\|_{\dot H^{t}_p(\crnp)}.\tag5.14
$$

$2^\circ$  When $\Omega $ is a bounded $C^{1+\tau }$-domain, $\tau
\ge 1$, there holds for  $s,t<\tau $: 
$$
\|d^{s}u\|_{\dot H_p^{t}(\comega)}\le C \|u\|_{\dot H^{t}_p(\comega)}.\tag5.15
$$
\endproclaim

\demo{Proof} $1^\circ$.
For integer $s$, the property is obvious since $x_n^s$ is then a
polynomial, a $C^\infty $-function on $\rn$. A noninteger power can
be written as $x_n^kx_n^s$ with $s\in (0,1)$ and $k\in \N_0$. When the  property has
been proved for $s\in (0,1)$, it follows immediately for $x_n^kx_n^s$
since $x_n^k$ is polynomial. It remains to treat the case $s\in
(0,1)$.

Since $x_n^s$ is a bounded function  on
$\ol B$, (5.14) holds for $t=0$. For $t=1$, consider $u\in
\dot H_p^1(\crnp)$. Here
$$
\partial_n(x_n^su)=x_n^s\partial_nu+sx_n^{s-1}u.\tag5.16
$$
For the first term,
$$
\|x_n^{s}\partial_nu\|_{L_p(B)}\le \|x_n^{s}\|_{L_\infty (B )} \|u\|_{\dot H_p^{1}(\crnp)} \le C \|u\|_{\dot H_p^{1}(\crnp)}. 
$$
For the second term,
by Theorem 5.4, since $s-1<0$, $s>0$,
$$
\|x_n^{s-1}u\|_{L_p(B )}\le C_1 \|u\|_{\dot H_p^{1-s}(\crnp)} \le
C_2 \|u\|_{\dot H_p^{1}(\crnp)} ,
$$
so altogether,
$$
\|\partial_n(x_n^{s}u)\|_{L_p(B )}\le C_3 \|u\|_{\dot H_p^{1}(\crnp)}. 
$$
With similar (easier) estimates of $\partial_j(x_n^su)$, $j<n$, we
conclude that $x_n^su\in \dot H_p^1(\crnp)$.

When the property is known for $t=k$ up to some $k\in\N$, we get it for
$t=k+1$ by using (5.16) and applying the result for $t=k$:
Let  $u\in \dot H_p^{k+1}(\crnp)$. Then
$$
\|x_n^s\partial_nu\|_{ \dot
H_p^{k}(\crnp)}\le C
\|\partial_nu\|_{ \dot
H_p^{k}(\crnp)}\le C' \|u\|_{ \dot
H_p^{k+1}(\crnp)},
$$
 by the result known for $\dot H_p^{k}(\crnp)$. And by Theorem 5.4,
$$
\|x_n^{s-1}u\|_{\dot H_p^k(\crnp)}\le C_1 \|u\|_{\dot
H_p^{k+1-s}(\crnp)}\le C_2\|u\|_{\dot
H_p^{k+1}(\crnp)}.
$$
It follows that
$$
\|\partial_n(x_n^{s}u)\|_{\dot H_p^k(\crnp)}\le C_3\|u\|_{\dot
H_p^{k+1}(\crnp)}.
$$
There are similar estimates of the derivatives $\partial_j$ for $j<n$,
so we
conclude that
$$
\|x_n^{s}u\|_{\dot H_p^{k+1}(\crnp)}\le C_4
\|u\|_{\dot H_p^{k+1}(\crnp)}.
$$

Thus the property (5.14) holds for arbitrarily high integer $t$.
Finally, the continuity of
the mapping $x_n^s$:
$$
u\in \{w\in\dot H_p^{t}(\crnp)\mid \supp w\subset \ol B \} \mapsto
x_n^{s}u\in \{w\in\dot H_p^{t}(\crnp)\mid \supp w\subset \ol B \} \tag5.17
$$
 for integer $t\ge 0$, extends to noninteger  $t>0$ by interpolation.

$2^\circ$. The statement carries over to curved sets $\Omega $ in the
same way as in the proof of Theorem 5.4. \qed
\enddemo

It is of interest to observe that the analysis in  \cite{T01, Prop.\
5.7} in fact covers a large family of spaces, namely the spaces
$F_{pq}^t(\rn)$ for $0< p\le \infty $, $0<q\le\infty $ with
$q=\infty $ if $p=\infty $, $t>n(\frac1p-1)_+$. Here is it shown that
when $\supp u\subset \ol B$, then
$$
\|x_n^{-t}u\|_{L_p(\rnp)}\le C\|f\|_{\dot F_{pq}^t(\crnp)}.\tag5.18
$$
For some of these spaces, we can extend the argumentation in
Proposition 5.3 and Theorem 5.4 to obtain generalizations of (5.5).
We use that  in these scales, $\ang D=\OP(\ang \xi )$
 defines a homeomorphism of $F_{pq}^{t+1}(\rn)$ onto
 $F_{pq}^{t}(\rn)$; so  $u\in \dot F_{pq}^{t+1}(\crnp)$ if and only if
 $u,\partial_1u,\dots,\partial_nu$ are in  $\dot
 F_{pq}^{t}(\crnp)$.

\proclaim{Theorem 5.6} Let $p\in (1,\infty )$ and $q\ge 2$, or
$p=q=\infty $.
Let $s,t\in\rp$.

$1^\circ$ There holds for functions $u$ supported in $ \ol B$:
$$
\|x_n^{-s}u\|_{\dot F_{pq}^{t}(\crnp)}\le C \|u\|_{\dot
F_{pq}^{s+t}(\crnp)}.
\tag5.19
$$

$2^\circ$ There holds on a bounded $C^\infty $-domain $\Omega $:
$$
\|d^{-s}u\|_{\dot F_{pq}^{t}(\comega)}\le C \|u\|_{\dot
F_{pq}^{s+t}(\comega)}.
\tag5.20
$$

\endproclaim

\demo{Proof} 
$1^\circ$. When $p\in (1,\infty )$ and  $q\ge 2$,
$H^t_p(\rn)=F^t_{p2}(\rn)$ is continuously embedded in $ F^t_{pq}(\rn)$, with
similar  embeddings for spaces over domains; in particular,
$$
L_p(\rn)=F^0_{p2}(\rn)\subset 
F^0_{pq}(\rn),\quad
\|u\|_{F^0_{pq}(\crnp)}\le C\|u\|_{L_p(\crnp)}.
$$
Then (5.18)  has the consequence
$$
\|x_n^{-s}u\|_{F^0_{pq}(\crnp)}\le C\|u\|_{F^s_{pq}(\crnp)}, \text{
when }\supp u\subset \ol B, \;s\ge 0.\tag5.21
$$
Departing from this estimate, the proofs given in Proposition 5.3 and
 Theorem 5.4, lifting to higher-order spaces by iteration (followed by
 interpolation), go over verbatim to the $F_{pq}^t$-scale, by replacing $H_p^t$ by
 $F_{pq}^t$. Let us just make the first nontrivial step explicit: Here we must
 show:
$$
u\in \dot F_{pq}^2(\crnp)\implies x_n^{-1} u\in  \dot F_{pq}^1(\crnp).\tag5.22
$$
Let $u\in \dot F_{pq}^2(\crnp)$. By (5.21),
$x_n^{-2}u\in \dot F_{pq}^0(\crnp)$. Note also that since $u\in
\dot F_{pq}^1(\crnp)$, 
$x_n^{-1}u\in \dot F_{pq}^0(\crnp)$ by (5.21).
 Moreover, since
$\partial_nu\in \dot F_{pq}^1(\crnp)$,
$x_n^{-1}\partial_nu\in \dot F_{pq}^0(\rnp)$ by (5.21). Then also 
$$
\partial_n(x_n^{-1}u)=-x_n^{-2}u+x_n^{-1}\partial_nu\in \dot F_{pq}^0(\crnp).
$$
For $j<n$, we simply have $
\partial_j(x_n^{-1}u)=x_n^{-1}\partial_ju\in \dot F_{pq}^0(\crnp)$,
 since $\partial_ju\in  \dot F_{pq}^1(\crnp)$. Altogether, since
 $x_n^{-1}u$ and its first derivatives are in $\dot F_{pq}^0(\crnp)$,
 (5.22) is obtained.

When $p=q=\infty $, i.e., $F_{pq}^t(\rn)=C_*^t(\rn)$, it is known that
$$
\{w\in e^+L_\infty (\rnp)\mid \supp w\subset \ol B\}\subset \{w\in \dot C^0_*(\crnp)\mid \supp w\subset \ol B\},
$$
so (5.18) implies that when $\supp
u\subset \ol B$,
$$
\|x_n^{-s}u\|_{\dot C_*^0(\crnp)}\le C\|u\|_{\dot C_*^s(\crnp)}.
$$
The method of Proposition 5.3 and Theorem 5.4  here leads to the estimates for function $u$ with $\supp u\subset \ol B$:
$$
\|x_n^{-s}u\|_{\dot C_*^{t}(\crnp)}\le C \|u\|_{\dot
C_*^{s+t}(\crnp)},
\tag5.23
$$
the version of (5.19) for $p=q=\infty $.

$2^\circ$. It is accounted for in \cite{T01} how the spaces over
 smooth bounded sets
 $\Omega $ can be defined from the half-space case by localization; hereby the estimates are
 transferred to these spaces.
 \qed 
\enddemo

\example{Remark 5.7}
The result for H\"older-Zygmund spaces $C^t_*$ sheds more light on the
embeddings considered in
Proposition 5.1. For sufficiently smooth domains it fills out the lacking integer-order cases: When $u$ is
supported in $\ol B$, then for $s,t\ge 0$,
$$
u\in \dot C_*^{s+t}(\crnp)\implies u\in x_n^s\dot C_*^{t}(\crnp).\tag5.24 
$$
This carries over to bounded domains by the usual localization tools,
showing that when $\Omega $ is bounded $C^{1+\tau }$, $\tau \ge 1$,
$$
\dot C_*^{s+t}(\comega)\subset d^s\dot
C_*^{t}(\comega)\text{ for }s,t\ge 0,\, s+t<\tau .\tag5.25 
$$

Note moreover that since $e^+x_n^s$ is in $C^s(\ol B)$, there holds
when $\supp u\subset \ol B$, $0<\varepsilon\le s<1$,
$$
u\in \dot C_*^{\varepsilon}(\crnp)\implies x_n^su\in \dot C_*^{\varepsilon}(\crnp).\tag5.26 
$$
This property can be lifted (using (5.23)) to all spaces
$C_*^t(\crnp)$ with $t\ge \varepsilon$ by the method of proof of
Theorem 5.6. E.g., the first step is:

We lift from $\dot C_*^\varepsilon(\crnp)$ to $\dot
C_*^{1+\varepsilon}(\crnp)$ by using (5.16), where the first term is
in $\dot C_*^\varepsilon(\crnp)$ since $\partial_nu$ is there, and the second
term is estimated using (5.23) by
$$
\|x_n^{s-1}u\|_{\dot C^\varepsilon(\crnp )}\le C_1\|u\|_{\dot
C^{\varepsilon+1-s}(\crnp )}\le C_2 \|u\|_{\dot
C^{1+\varepsilon}(\crnp )}.
$$
Thus $\partial_n(x_n^su)\in \dot C_*^{\varepsilon}(\crnp)$, and with
similar estimates for $\partial_j(x_n^su)$ and $x_n^su$, we conclude
that $x_n^{s}u\in \dot C^{1+\varepsilon}(\crnp )$.

The complete result is that for $u$ with $\supp u\subset \ol B$:
$$
\|x_n^su\|_{\dot C_*^{t}(\crnp)}\le C \|u\|_{\dot
C_*^{t}(\crnp)},\text{ when }s,t\in \rp.\tag5.27
$$
The result extends by localization to spaces $\dot C_*^t(\comega)$ with
$\Omega$  bounded and $C^{1+\tau }$, $\tau \ge 1$:
$$
\|d^su\|_{\dot C_*^{t}(\comega)}\le C \|u\|_{\dot
C_*^{t}(\comega)},\text{ when }s,t\in (0,\tau ).\tag5.28
$$
\endexample

\subhead 6.  Gradient estimates
\endsubhead

An interesting question 
is the
validity of weighted gradient estimates, both of
$\nabla u$ and  $\nabla (u/d^a)$ multiplied by powers of $d$, when $u$ is a solution of
(1.1). The results in the following are new for solutions of (1.1)
with $x$-dependent operators $P$; it is moreover new that they can be
obtained in Bessel-potential spaces.

Fall and Jarohs showed in \cite{FJ21} some gradient estimates for
solutions of (1.1) with $P=(-\Delta )^a$ in the framework of
H\"older spaces: When $a> \frac12$ and $\Omega $ is bounded $C^{1,1}$, then
when $f\in L^\infty (\Omega )$, $\beta <2a-1$,
$$
\|d^{1-a}\nabla u\|_{C^\beta (\comega)^n }\le C (\|f\|_{L_\infty (\Omega
)}+ \|u\|_{L_\infty (\Omega )}).\tag6.1
$$
Furthermore, for any $\varepsilon >0$,
$$
\|d^{1-a+\varepsilon }\nabla (u/d^a)\|_{L_\infty (\Omega)^n }\le C \|f\|_{L_\infty (\Omega )}+\|u\|_{L_\infty (\Omega )}).\tag6.2
$$
Related estimates were found earlier in
Ros-Oton and Serra \cite{RS14, Th.\ 1.4}
and in Ros-Oton, Serra and Valdinoci  \cite{RSV17, Prop.\ 5.1}, the
latter including  a larger class of $x$-independent
 operators $P$. 
Estimates in Bessel-potential spaces with some of the same flavour, but with powers $(-\Delta
)^{t}$ instead of $\nabla$, were worked out by  Abdellaoui et al.\ in
\cite{AFLY22}. Interior gradient estimates have  been
shown by Kuusi, Nowak and Sire \cite{KNS25} in spaces of Sobolev
and Besov
type, for a class of
fractional-order operators with some overlap with the one we consider,
but with a different aim: to allow measures $f$.

The estimates (6.1), (6.2) imply that solutions $u$ of (1.1) with
$f\in L_\infty (\Omega )$ satisfy:
$$
d^{1-a}\nabla u\in C^\beta(\comega)^n \text{ with } d(x)^{1-a}\nabla
u(x)\cdot \nabla d(x) = a u(x)/d(x)^a \text{ on }\partial\Omega .\tag6.3 
$$

We shall now show how such estimates can be generalized, and also
applied to $x$-dependent operators $P$ (as in Hypothesis 2.1), both in the
H\"older-Zygmund scale and the Bessel-potential scale of spaces.

It turns out that (6.1) and its generalizations  are (after our preparations with embedding rules in
Section 5) the easiest to come
by, since they only draw on Poisson operators in the standard
Boutet de Monvel calculus, whereas (6.2) and its generalizations need
some further refined calculations. So let us begin with (6.1). We can
here include higher-order versions with an extra smoothness parameter $t$.

\subsubhead 6.1 Estimates of $\nabla u$ in H\"older-Zygmund spaces \endsubsubhead

\proclaim{Theorem 6.1}  Let $a\in (\frac12,1)$, let $t\in [0,1-a)$, and let $\Omega $ be bounded
 $C^\infty $.
 Then
 there is a constant $C>0$ such that for $u\in
 C_*^{a(2a+t)}(\comega)$,
 $$
\|d^{1-a}\nabla u\|_{\ol C_*^{2a-1+t} (\Omega )^n}\le C \|u\|_{C_*^{a(2a+t)}(\comega)}.\tag6.4
$$

\endproclaim

\demo{Proof}
Decompose $u$  in $u=w+z$ as in 
Corollary 3.8 with $s=t\in [0, 1-a)$; here 
$$
w\in \dot C_*^{2a+t}(\comega),\quad z=d^aK_{(0)}\varphi \text{ for a
}\varphi \in C_*^{a+t}(\partial\Omega),  
$$
 with norms bounded in terms of $\|u\|_{C_*^{a(2a+t)}(\comega)}$. 

For the function $w$ we have 
 since $2a-1>0$ and $1-a>0$, using (5.28),
$$
\|d^{1-a}\nabla w\|_{\dot C_*^{2a-1+t}(\comega)^n}\le C_1\|\nabla
w\|_{\dot C_*^{2a-1+t}(\comega)^n}\le C_2\|w\|_{\dot C_*^{2a+t}(\comega)}\le
C_3\|u\|_{C_*^{a(2a+t)}(\comega)}.\tag6.5
$$

For the function $z=d^aK_{(0)}\varphi $ we shall use that $K'=d\nabla
K_{(0)}$ is a Poisson operator of order 0 in the Boutet de Monvel
calculus. In fact, in local
coordinates reduced to the $\rnp$-situation, $d$ equals $x_n$, and the symbol estimates for
$x_n\partial_jK_{0}$ are those of a Poisson operator of order 0. (The
basic calculation on $\rnp$ is recalled below in Theorem 6.9, proof of
(6.24).)
Now 
$$
d^{1-a}\nabla(d^aK_{(0)}\varphi )=d\nabla K_{(0)}\varphi
+d^{1-a}(\nabla d^a)K_{(0)}\varphi =K'\varphi
+a(\nabla d)K_{(0)}\varphi\in \ol C_*^{a+t}(\Omega ),\tag6.6 
$$
with norm bounded in terms of $\|\varphi \|_{C_*^{a+t}(\partial\Omega )}$, since $K'$
and $K_{(0)}$ are Poisson operators of order 0. It follows that
$$
\|d^{1-a}\nabla z\|_{\ol C_*^{a+t}(\Omega )^n}
\le C_4\| \varphi \|_{ C_*^{a+t}(\partial\Omega )}\le C_5\|u\|_{C_*^{a(2a+t)}(\comega)}.
$$

The sum $u=w+z$ then satisfies (6.4), since $a> 2a-1$.
\qed
\enddemo

For higher $t$, Theorem 3.9 can be used in a similar way to
prove the inequality (6.4):

\proclaim{Theorem 6.2}  Let $a\in (\frac12,1)$,  and let $\Omega $ be bounded
 $C^\infty $. The estimate {\rm (6.4)} holds for all  $t\ge 0$.
\endproclaim

\demo{Proof} Let $M$ be a large positive integer, and let $t\in
(-a+M-1,- a+M)$. Then by Theorem 3.9, any $u\in C_*^{a(2a+t)}$ has a unique decomposition 
$$
u=w+z,\; w\in \dot C_*^{2a+t}(\comega), \; z=d^a\widetilde {\Cal K}'_M\varphi 
, $$
where $\varphi =\varrho _M^au=\{\gamma _0^au, \dots, \gamma _{M-1}^au\}\in \prod
_{j<M}C_*^{t+a -j}(\partial\Omega )$, and 
$\widetilde {\Cal K}'_M$ is a vector of standard Poisson operators of
orders $\{0,-1,\dots,- M+1\}$. For the function $w$ we have the
estimate (6.5). For $z$, we have: 
$$
\aligned
d^{1-a}\nabla z&=d^{1-a}\nabla(d^a\widetilde{K}'_M\varphi
)=d\nabla\widetilde{K}'_M\varphi
+d^{1-a}\nabla(d^a)\widetilde{K}'_M\varphi \\
&=d\nabla \widetilde
{\Cal K}'_M
\varphi 
+a(\nabla d)\widetilde {\Cal K}'_M\varphi , 
\endaligned
\tag6.7
$$
where each component of $d\nabla \widetilde
{\Cal K}'_M$ is a vector of Poisson operators of the same orders (and with the same continuity
properties) as $\widetilde
{\Cal K}'_M$.
It follows that
$$
\|d^{1-a}\nabla z\|_{\ol C_*^{a+t}(\Omega )^n}
\le C_4\| \varphi \|_{\prod _{j<M}C_*^{a+t -j}(\partial\Omega ) }\le C_5\|u\|_{C_*^{a(2a+t)}(\comega)}.
$$

For $u=w+z$ we then find (6.4) for the considered $t$, since $a>2a-1$. The continuity
of the mapping
$$
u\in C_*^{a(2a+t)}(\comega)\mapsto d^{1-a}\nabla u\in \ol C_*^{2a-1+t}(\Omega )^n
$$
has then been obtained both for an arbitrarily high $t$ and for $t$
close to 0
(by Theorem 6.1); then it follows for general $t>0$ by interpolation. \qed
\enddemo

This leads to the result for solutions of (1.1):

\proclaim{Corollary 6.3}  Let $a\in (\frac12,1)$, let $t\ge 0$, let $P$
 satisfy Hypothesis {\rm 2.1} with $\tau =\infty $,
 and let $\Omega $ be bounded
 $C^\infty $. When $u$ is a solution of {\rm (1.1)} with $f\in \ol
 C_*^t(\Omega )$, then
 $$
 \aligned
\|d^{1-a}\nabla u\|_{\ol C_*^{2a-1+t} (\Omega )^n}&\le C (\|f\|_{\ol
C_*^t(\Omega )}+\|u\|_{L_\infty (\Omega )})\\
\|d^{1-a}\nabla u\|_{\ol C_*^{2a-1} (\Omega )^n}&\le C (\|f\|_{L_\infty (\Omega )}+\|u\|_{L_\infty (\Omega )})
\endaligned
\tag6.8
$$
When there is uniqueness of solution, the terms  $\|u\|_{L_\infty
(\Omega)}$ can be omitted.
\endproclaim

\demo{Proof}  By the theory in \cite{G14} (cf.\ Th.\ 3.2 and Example
3.3 there), (1.1) is Fredholm solvable for   $f\in \ol
C_*^t(\Omega )$, $u\in \dot C_*^\sigma $
(some $\sigma \in (a-1,a]$), with solutions $u$ in $C_*^{a(2a+t)}(\comega)$; they satisfy
$$
\|u\|_{C_*^{a(2a+t)}(\comega)}\le C_1(\|f\|_{\ol
C_*^t(\Omega )}+\|u\|_{\dot
C_*^\sigma (\comega )}).\tag6.9
$$
When there is uniqueness of solution,
$\|u\|_{C_*^{a(2a+t)}(\comega)}\simeq \|f\|_{\ol
C_*^t(\Omega )}$.

For $t=0$, we use  $L_\infty (\Omega )\subset\ol
C_*^0(\Omega )$,  together with $L_\infty (\Omega )\subset\dot
C_*^\sigma (\comega )$ (when $\sigma \le 0$), to see that
$$
\|u\|_{C_*^{a(2a)}(\comega)}\le C_2(\|f\|_{L_\infty (\Omega )}+\|u\|_{L_\infty (\Omega )}).\tag6.10
$$
Then  the estimate in Theorem 6.1 implies (6.8) for $t=0$. It
follows similarly for $t>0$ using Theorem 6.2. \qed

\enddemo

In the  case of smooth domains this sharpens the estimate (6.1) shown in
\cite{FJ21} for  $(-\Delta )^a$, where $\beta <2a-1$.
After we told this to the authors,
they sent a personal
communication where they can extend their proof to cover $\beta =2a-1$;
they also correct an assumption (1.8) in their paper.

The estimates in $a$-transmission spaces carry over to nonsmooth domains by localization:

\proclaim{Theorem 6.4}  Let $a\in (\frac12,1)$, $\tau \ge 1$,  $t\ge 0$,
 and let $\Omega $ be 
 a bounded
 $C^{1+\tau } $-domain. Assume that 
$$
\aligned
\text{ if }t\le 1-a, \text{ then }\tau >t+2a,\\
\text{ if }t>1-a, \text{ then }\tau >t+2a+2.
\endaligned\tag6.11
$$
 Then
$$
\|d^{1-a}\nabla u\|_{\ol C_*^{2a-1+t} (\Omega  )^n}\le C
\|u\|_{C_*^{a(2a+t)}(\comega)}.
\tag6.12
$$
\endproclaim

\demo{Proof} The result in the smooth case, Theorem 6.2, holds in
particular for smooth subsets $\Omega _R$ of $\rnp$, where
$B'_R=\{x\mid |x'|<R, x_n=0\}$ is part of the boundary $\partial\Omega
_R$, and
$\rnp\cap B_R\subset \Omega _R\subset \rnp$. Here we have
the inequality (6.12) in particular when $u\in C_*^{a(2a+t)}(\comega_R)$ is supported in
$\ol B_{R/2}\cap \crnp$, say.

For a given bounded $C^{1+\tau }$-domain $\Omega $ we use localization
as described in Proposition 4.3 to show that the inequality carries
over to a coordinate patch  $U\cap \R^n_\zeta $ using the diffeomorphism $F_\zeta
$, when the conditions in (6.11) hold, cf.\ Corollary A.3. The information is pieced
together from finitely many such patches by cut-off functions and
rotations.
It is used that the differentiations $\partial_j$
are turned into first-order differentiations with $C^\tau
$-coefficients by the diffeomorphisms; multiplication by these
coefficients are bounded
in the considered norms.
\qed
\enddemo

The results that we have for the Dirichlet problem (1.1) in $C_*^s$-spaces over
nonsmooth domains, quoted in Theorem 2.3 $2^\circ$, carry a ``loss of
$\varepsilon $'' since they are derived from estimates in Bessel-potential spaces by use of Sobolev embeddings.
However, the
$H^s_q$-results have a small advantage in that the regularity
parameter $t$ need only satisfy $t<\tau -2a$, not the stronger
condition (6.11). To use this, we shall base the estimates directly on
$H^s_q$-estimates. The final result will be given at the end, in
Theorem 6.19.

\subsubhead 6.2 Estimates of $\nabla u$ in Bessel-potential spaces \endsubsubhead

The same type of gradient estimates in Bessel-potential spaces can be shown by
very similar methods, thanks to the embedding rules worked out in Section 5.

Recall that $q\in (1,\infty )$ throughout. We shall show gradient estimates for functions in the solution space
$H_q^{a(2a+t)}(\comega)$, $\Omega $ bounded or of half-space type.

First there is the result in the case of $C^\infty $-domains.

\proclaim{Theorem 6.5}  Let $a\in (\frac12,1)$,  and let $\Omega $ be bounded
 $C^\infty $. There is a $t_0$ with  $1-2a\le t_0<0$ such that for all $t\ge t_0$, there holds
 $$
\|d^{1-a}\nabla u\|_{\ol H_q^{2a-1+t} (\Omega )^n}\le C \|u\|_{H_q^{a(2a+t)}(\comega)}.\tag6.13
$$

\endproclaim

\demo{Proof} Note first that for $w$ in the space $\dot
H_q^{2a+t}(\comega)$, there holds when $t\ge 1-2a$, by use of the
rule (5.15):
$$
\|d^{1-a}\nabla w\|_{\dot H_q^{2a-1+t}(\comega)^n}\le C_1\|\nabla
w\|_{\dot H_q^{2a-1+t}(\comega)^n}\le C_2\|w\|_{\dot H_q^{2a+t}(\comega)}.
\tag6.14
$$

The spaces $H_q^{a(2a+t)}(\comega)$ are defined for
all $t>-a-\frac1{q'}$.
For $t$ in the interval $I_0=(-a-\frac1{q'}, -a+\frac1q)$ (of length 1), 
$H_q^{a(2a+t)}(\comega)=\dot H_q^{2a+t}(\comega)$, cf.\ (2.15).

If $a\le \frac1q$, the right endpoint of the interval $I_0$ is $\ge 0$ (but $<1$). Then
there is a $t_0<0$ in the interval satisfying $t_0\ge 1-2a$. For this
$t_0$, we have from (6.14) since $H_q^{a(2a+t_0)}(\comega)=\dot
H_q^{2a+t_0}(\comega)$ that
$$
\|d^{1-a}\nabla u\|_{\dot H_q^{2a-1+t_0}(\comega)^n}
\le C\|u\|_{H_q^{a(2a+t_0)}(\comega)}.\tag6.15
$$

If $a>\frac1q$, the next interval $I_1=(-a+\frac1{q}, -a+1+\frac1q)$
contains 0 and hence contains a number $t_0<0$ satisfying $t_0\ge 1-2a$.
For this value $t_0$, we apply Corollary 3.3 with $s=t_0$. Decompose
$u$  in $u=w+z$; here 
$$
w\in \dot H_q^{2a+t_0}(\comega),\quad z=d^aK_{(0)}\varphi \text{ for a
}\varphi \in B_q^{a+t_0-1/q}(\partial\Omega),  
$$
 with norms bounded in terms of $\|u\|_{H_q^{a(2a+t_0)}(\comega)}$. For
 $w$ we have the estimate (6.14), which implies
$$
\|d^{1-a}\nabla w\|_{\dot H_q^{2a-1+t_0}(\comega)^n}\le
C\|u\|_{H_q^{a(2a+t_0)}(\comega)}.\tag6.16
$$
For $z=d^aK_{(0)}\varphi$ we use the observation in the proof of
Theorem 6.1, that 
$$
d^{1-a}\nabla(d^aK_{(0)}\varphi )=K'\varphi
+a(\nabla d)K_{(0)}\varphi,
$$
where $K'$ and $K_{(0)}$ are standard Poisson operators of order 0, so that 
$$
\|d^{1-a}\nabla z\|_{\ol H_q^{a+t_0}(\Omega )^n}
\le C\| \varphi \|_{ B_q^{a+t_0-1/q}(\partial\Omega )}\le C'\|u\|_{H_q^{a(2a+t_0)}(\comega)}.
$$
Hence the sum $u=w+z$ satisfies (6.13) for $t=t_0$.

Now let $M$ be an integer $>1$, and let $t_1\in I_M=(-a+M-\frac1{q'},
-a+M+\frac1q)$. Here we have the decomposition from Theorem 3.5, for
$u\in H_q^{a(2a+t_1)}$: 
$$
u=w+z,\; w\in \dot H_q^{2a+t_1}(\comega), \; z=d^a\widetilde {\Cal K}'_M\varphi 
, $$
where $\varphi =\varrho _M^au=\{\gamma _0^au, \dots, \gamma _{M-1}^au\}\in \prod
_{j<M}B_q^{a+t_1 -j-1/q}(\partial\Omega )$, and $\widetilde {\Cal K}'_M$ is a vector of standard Poisson operators of
orders $\{0,-1,\dots,- M+1\}$. Here we can use the calculation (6.7) to
see that
$$
\|d^{1-a}\nabla z\|_{\ol H_q^{a+t_1}(\Omega )^n}
\le C\| \varphi \|_{\prod _{j<M}B_q^{a+t_1 -j-1/q}(\partial\Omega ) }\le C'\|u\|_{H_q^{a(2a+t_1)}(\comega)},
$$
and the desired estimate of $u$ follows by addition.

Now (6.13) has been shown for $t$ equal to $t_0$ and $t_1$, and then
it follows for $t\in [t_0,t_1]$ by interpolation. Since $M$, hence
$t_1$, can be taken arbitrarily large, the estimate holds for all
$t\in [t_0,\infty )$. \qed

\enddemo

Just like for $C_*$-spaces, we can draw consequences for nonsmooth domains
and for actual  solutions of (1.1).

\proclaim{Theorem 6.6} Let $a\in (\frac12,1)$, $\tau \ge 1$,  
 and let $\Omega $ be 
 a bounded
 $C^{1+\tau } $-domain. Let $t_0\in [1-2a,0)$ be chosen as in Theorem
 {\rm 6.5}. For $t\in [t_0,\tau -2a)$ there holds
$$
\|d^{1-a}\nabla u\|_{\ol H_q^{2a-1+t} (\Omega  )^n}\le C
\|u\|_{H_q^{a(2a+t)}(\comega)}.
\tag6.17
$$
\endproclaim

\demo{Proof} The proof goes as in the proof of Theorem
6.4, using that the localization of the spaces is allowed when $t\in
(-a-\frac1{q'}, \tau -2a)$. \qed
\enddemo

\proclaim{Corollary 6.7} Let $a\in (\frac12,1)$, let $P$
 satisfy Hypothesis {\rm 2.1} with $\tau \ge 1 $,
 and let $\Omega $ be bounded
 $C^{1+\tau } $. Let $t\in [0,\tau -2a)$.

When $u$ is a solution of {\rm (1.1)} with $f\in \ol
 H_q^t(\Omega )$, then
 $$
\|d^{1-a}\nabla u\|_{\ol H_q^{2a-1+t } (\Omega )^n}\le C (\|f\|_{\ol
H_q^t(\Omega )}+\|u\|_{L_q(\Omega )}).
\tag6.18
$$
When there is uniqueness of solution, the term  $\|u\|_{L_q
(\Omega)}$ can be omitted.
\endproclaim

\demo{Proof} Here we refer to the information in Theorem 2.2, showing that when $u$ solves (1.1) with $f\in \ol
H_q^t(\Omega )$ for a $t\ge 0$, then
$\|u\|_{H_q^{a(2a+t )}(\comega)}$ is finite, equivalent to $\|f\|_{\ol
H_q^t(\Omega )}$ in the case with uniqueness, and generally dominated by $\|f\|_{\ol
H_q^t(\Omega )}+\|u\|_{L_q (\Omega )}$. This is combined with the estimate
in Theorem 6.6. \qed
\enddemo

\subsubhead 6.3 Estimates of $\nabla( u/d^a)$ in Bessel-potential spaces \endsubsubhead

Estimates of the type (6.2) are a little more demanding, since the
 occurring Poisson-type operators are not standard operators in the
 Boutet de Monvel calculus.

We recall that $H_q^{a(2a+t)}(\crnp)$ equals $\dot H_q^{2a+t}(\crnp)$
for $t\in (-a-\frac1{q'},-a+ \frac1q)$, which overlaps with $t\ge 0$ when
$a<\frac1q$. When $t\in (-a+\frac1q,-a+ 1+\frac1q)$, the space has the decomposition in Corollary 3.3,
where there is a component $z=d^aK_{(0)}\varphi
$, $\varphi =\gamma _0(u/d^a)$, besides the component in $\dot
H_q^{2a+t}(\crnp)$.

In this paper, we shall not work out the details for higher
values of $t$; they can be based on Theorem 3.5.

The part $w\in \dot H_q^{2a+t}(\crnp)$ can be estimated by use of
the embedding rules in Section 5. For the part $z$, we need the
following analysis.

 First recall some precise formulas for
 Fourier transforms:
Denoting $1|_{\rp}=H(x_n)$ (the Heaviside function), we have for
$\sigma >0$, $t> -1$, 
$$
\align
{\Cal F}_{x_n\to \xi _n }(H(x_n)e^{-\sigma x_n})&=\frac{1}{\sigma
  +i\xi _n },\tag6.19 \\
    {\Cal F}_{x_n\to \xi _n }(H(x_n)x_n^te^{-\sigma x_n})&=\frac{c}{(\sigma
  +i\xi _n )^{t+1}},\quad c=\Gamma
  (t+1).\tag6.20
\endalign
$$
  The complex number $\sigma +i\xi _n $ has real part $\sigma >0$, so its
noninteger powers (defined with a cut along the negative axis $\rmi$)
make good sense. Formula (6.19) is elementary; proofs of formula
(6.20) are found e.g.\ in Schwartz \cite{S61}, (V,1;44), and in the lines
after Example 7.1.17 in \cite{H83} (with different conventions).

\example{Remark 6.8} In the model case where $\Omega =\rnp$ and $P=(1-\Delta )^a$ with  $a\in (0,1)$,  (1.1) 
has when $t>-a-\frac1{q'} $ a unique solution $u\in H_q^{a(2a+t)}(\crnp)$
for any $f\in \ol H_q^t(\rnp)$, so $r^+P\colon
H_q^{a(2a+t)}(\crnp)\simto \ol H_q^t(\rnp)$ is a
homeomorphism.

Consider it in particular for $a\in (\frac12,1)$ and $t+a \in(\frac1q,
1+\frac1q)$. Here there is by Lemma {\rm 3.1} with $\mu =a$ a decomposition:
$$
u=w+z,\quad w\in \dot
H_q^{2a+t}(\crnp), \quad z= K^a _0 \varphi \text{ with }\varphi  \in B_q^{t+a -1/q}(\R^{n-1} )  ,\tag 6.21
$$
where we recall that 
$$
 K^a  _0= \tfrac1{\Gamma (a +1)}x_n^a  e^+ K_0.\tag6.22
$$
The norms $\|w\|_{\dot H_q^{2a+t}(\crnp)}$ and $\|\varphi \|_{B_q^{t+a
-1/q}(\R^{n-1} ) }$ are bounded in terms of
$\|u\|_{H_q^{a(2a+t)}(\crnp)}$ (and equivalently in terms of  $\|f\|_{\ol H_q^t(\rnp)}$).
\endexample

Then we can show gradient estimates for  $z/x_n^a$:

\proclaim{Theorem 6.9} Let $a\in (0,1)$. For $z=K_0^a\varphi $
on $\rnp$, there holds 
when $\frac1{q}< a+t<1+\frac1q$, $s\ge 0$,
$$
\|x_n^{1-a+s}\nabla (z/x_n^a)\|_{\dot H_q^{s+t}(\crnp)^n}\le  C\|\varphi
\|_{B_q^{a+t-1/q}(\R^{n-1})}.\tag 6.23
$$

There also holds, for all $t\in\R$:
$$
\|x_n^{1-a}\nabla z\|_{\ol H_q^{a+t}(\rnp)^n}\le C\|\varphi
\|_{B_q^{a+t-1/q}(\R^{n-1})}.\tag 6.24
$$
\endproclaim

\demo{Proof}
We shall use some precise facts around the Poisson operator
$K_0=\operatorname{OPK} ((\ang{\xi '}+i\xi _n)^{-1})$. Recall that when
$k(x',\xi )$ is a Poisson symbol, the corresponding operator is
defined by
$$
\operatorname{OPK}(k(x',\xi ))\psi =\Cal F^{-1}_{\xi \to x}[k(x',\xi )\hat \psi (\xi ')],
$$
which we can also write as 
$$
\widetilde{\operatorname{OPK}}(\tilde k(x',x_n,\xi '))\psi =\Cal
F^{-1}_{\xi '\to x'}[\tilde k(x',x_n,\xi ')\hat \psi (\xi ')],
$$
where $k(x',\xi )$ is the {\it symbol} and $\tilde k(x',x_n,\xi
')=\Cal F^{-1}_{\xi _n\to x_n}k(x',\xi )$ is the associated {\it
symbol-kernel}. In the case of $K_0$, the symbol is $(\ang{\xi
'}+i\xi _n)^{-1}$ and the symbol-kernel is $H(x_n)e^{-\ang{\xi '} x_n}$,
cf.\ (6.19). The notion of Poisson operators extends to some symbols
of the kind $(\ang{\xi
'}+i\xi _n)^{r}$ with noninteger $r$; let's call them Poisson-type
operators. In particular, 
$$
\aligned
K^a_0 &=\tfrac1{\Gamma (a+1)}x_n^aK_0
=\widetilde{\operatorname{OPK}}(\tfrac1{\Gamma
(a+1)}x_n^aH(x_n)e^{-\ang{\xi '} x_n}) =\operatorname{OPK}((\ang{\xi '}+i\xi _n)^{-a-1})\\
&=\operatorname{OP}((\ang{\xi '}+i\xi _n)^{-a})\operatorname{OPK}((\ang{\xi '}+i\xi _n)^{-1})
=\Xi _+^{-a}K_0 ,
\endaligned\tag6.25
$$
cf.\ (6.20). Note that here multiplication of the symbol-kernel by $\frac1{\Gamma
 (a+1)}x_n^a$ corresponds to multiplication of the symbol by $(\ang{\xi '}+i\xi
 _n)^{-a}$.

Begin with the estimate (6.24), which is obtained by tools within the
Boutet de Monvel calculus. In $\nabla z=(\partial_1z,\dots, \partial_nz)$, the most important
element is $\partial_nz$. For this we have:
$$
\aligned
x_n^{1-a}\partial_nz&=x_n^{1-a}\partial_n\Cal
F^{-1}_{\xi '\to x'}(\tfrac{1}{\Gamma
(a+1)}x_n^aH(x_n)e^{-\ang{\xi '}x_n}\hat\varphi (\xi '))\\
&= x_n^{1-a}\widetilde{\operatorname{OPK}}[\tfrac{1}{\Gamma
(a+1)}(ax_n^{a-1}-x_n^{a}\ang{\xi '})H(x_n)e^{-\ang{\xi '}x_n}]\varphi
\\
&= \tfrac a{\Gamma (a+1)}K_0\varphi -
\widetilde{\operatorname{OPK}}[\tfrac{1}{\Gamma
(a+1)}x_n\ang{\xi '}H(x_n)e^{-\ang{\xi '}x_n}]\varphi \\
&= \tfrac a{\Gamma (a+1)}K_0\varphi -
\operatorname{OPK}[\tfrac{\Gamma (2)}{\Gamma
(a+1)}\ang{\xi '}(\ang{\xi '}+i\xi _n)^{-2}]\varphi. \endaligned\tag6.26
$$
Here $K_0$ as well as the operator $\operatorname{OPK}[\ang{\xi
'}(\ang{\xi '}+i\xi _n)^{-2}]$ are standard Poisson operators of order 0, so
$\varphi \in B_q^{a+t -1/q}(\R^{n-1} )$ implies $x_n^{1-a}\partial_nz\in
\ol H_q^{a+t}(\crnp)$, for any $t\in\R$ (as known e.g.\ from \cite{G90}). The operators $x_n^{1-a}\partial_j $ ($j<n$) map in
a similar way, so altogether, 
$$
\|x_n^{a-1}\nabla z\|_{\ol H_q^{a+t}(\rnp)^n}\le C\|\varphi
\|_{B_q^{a+t-1/q}(\R^{n-1})}, \text{ any }t\in\R,
$$
showing (6.24).

For the estimate (6.23) we show for $s\ge 0$ by the tools recalled above,
$$
\aligned
x_n^{1-a+s}\partial_n(z/x_n^a)&=x_n^{1-a+s}\partial_n\tfrac{1}{\Gamma
(a+1)}K_0\varphi \\
&=x_n^{1-a+s}\partial_n\Cal
F^{-1}_{\xi '\to x'}
(\tfrac{1}{\Gamma
(a+1)}H(x_n)e^{-\ang{\xi '}x_n}\hat\varphi (\xi '))\\
&= x_n^{1-a+s}\widetilde{\operatorname{OPK}}[\tfrac{1}{\Gamma
(a+1)}(\delta (x_n)-\ang{\xi '})H(x_n)e^{-\ang{\xi '}x_n}]\varphi
\\
&=- \widetilde{\operatorname{OPK}}[\tfrac{1}{\Gamma
(a+1)}x_n^{1-a+s}\ang{\xi '}H(x_n)e^{-\ang{\xi '}x_n}]\varphi
\\
&=- 
\operatorname{OPK}[\tfrac{\Gamma (2-a+s)}{\Gamma
(a+1)}\ang{\xi '}(\ang{\xi '}+i\xi _n)^{-2+a-s}]\varphi\equiv K^{(n)}\varphi . \endaligned\tag6.27
$$
It is used that $x_n^{1-a+s}\delta (x_n)=0$. The Poisson-like operator
$K^{(n)}$ can be written:
$$
\aligned
K^{(n)}&=c_{a,s}\operatorname{OPK}[\ang{\xi '}(\ang{\xi '}+i\xi
_n)^{-1+a-s}(\ang{\xi '}+i\xi _n)^{-1}]\\
&=c_{a,s}\Xi _+^{a-1-s}K_0\ang{D'};\quad c_{a,s}=-\tfrac{\Gamma (2-a+s)}{\Gamma
(a+1)}.
\endaligned\tag6.28
$$

Here $\ang{D'}$ maps  $\varphi
\in B_q^{a+t-1/q}(\R^{n-1})$ to $B_q^{a+t-1-1/q}(\R^{n-1})$, which is mapped to
$\ol H_q^{a+t-1}(\rnp)$ by $K_0$. Now $\ol H_q^{a+t-1}(\rnp)$
identifies with $\dot H_q^{a+t-1}(\rnp)$ when $a+t\in (\frac1{q},1+\frac1q)$.
Then we can apply $\Xi _+^{a-1-s}$, which maps the space into $\dot
H_q^{s+t}(\crnp)$. Thus $K^{(n)}\varphi \in \dot H_q^{s+t}(\crnp)$.

To sum up,
$$
K^{(n)}\colon B_q^{a+t-1/q}(\R^{n-1})\to \dot H_q^{s+t}(\crnp),
\text{ when $a+t\in (\tfrac1{q},1+\tfrac1q)$, $s\ge 0$} .\tag6.29
$$ 

For $j<n$,
$$
\aligned
x_n^{1-a+s}\partial_j(z/x_n^a)&=x_n^{1-a+s}\tfrac{1}{\Gamma
(a+1)}\partial_jK_0\varphi \\
&=\widetilde{\operatorname{OPK}}
(\tfrac{x_n^{1-a+s}}{\Gamma
(a+1)}(i\xi _j)H(x_n)e^{-\ang{\xi '}x_n})\varphi \\
&= 
\operatorname{OPK}[\tfrac{\Gamma (2-a+s)}{\Gamma
(a+1)}(i\xi _j)(\ang{\xi '}+i\xi _n)^{-2+a-s}]\varphi\equiv K^{(j)}\varphi . \endaligned
$$
We write $K^{(j)}$ as 
$$
K^{(j)}
=-c_{a,s}\Xi _+^{a-1-s}K_0\partial_j,
\tag6.30
$$
and it is seen as in the treatment of $K^{(n)}$ that
$$
K^{(j)}\colon B_q^{a+t-1/q}(\R^{n-1})\to \dot H_q^{s+t}(\crnp)\text{
when $a+t\in (\tfrac1{q},1+\tfrac1q)$, $s\ge 0$}.\tag6.31 
$$

Collecting the terms, we find (6.23).\qed

\enddemo

Next, we show estimates of $\nabla(w/x_n^{a})$ when $w$ is in a
supported space $\dot H_q^{2a+t}(\crnp)$.

\proclaim{Lemma 6.10} Let $a\in (\frac12,1)$, $s \in [0,2a-1)$, and $t\ge 
 -s$ (in particular, $t>1-2a$).
When $w\in \dot H_q^{2a+t}(\crnp)$ is supported in $\ol B $,
there holds:
$$
\|x_n^{1-a+s }\nabla (w/x_n^a)\|_{\dot H_q^{s+t} (\crnp)^n}\le
C\|w\|_{\dot H_q^{2a+t}(\crnp)}. 
\tag6.32
$$

\endproclaim

\demo{Proof}
 Since $w\in \dot H_q^{2a+t}(\crnp)$, $\nabla w\in \dot
H_q^{2a-1+t}(\crnp)^n$, and hence by Theorem 5.4,
$$
\nabla w\in \dot H_q^{2a-1+t}(\crnp)^n\subset x_n^{2a-1-s }\dot
H_q^{s +t}(\crnp)^n,  
$$
with norm dominated by that of $w$ in $\dot H_q^{2a+t}(\crnp)$.
Then for $w/x_n^a$,
$$
x_n^{1-a+s }\nabla (w/x_n^{a})=
x_n^{1-2a+s }\nabla w -ax_n^{-2a+s }(0,\dots,0,w)
\in \dot H_q^{s+t} (\crnp)^n,
$$
by Theorem 5.4. \qed
\enddemo

This will now be combined with the estimates in Theorem 6.9 to obtain
gradient estimates for  $a$-transmission spaces.

\proclaim{Theorem 6.11}  Let $a\in (\frac12,1)$, $s\in [0,2a-1)$, $t\in
 (-a-\frac1{q'},1-a+\frac1q)$ and $t\ge -s$.
For $u\in H_q^{a(2a+t)}(\crnp) $ with $\supp u\subset \ol B$, there holds:
$$
\|x_n^{1-a+s}\nabla (u/x_n^a)\|_{\dot H_q^{s+t} (\crnp )^n}\le C \|u\|_{H_q^{a(2a+t)}(\crnp)}.\tag6.33
$$

\endproclaim

\demo{Proof}
When $t\in
 (-a-\frac1{q'},-a+\frac1q)$, $H_q^{a(2a+t)}(\crnp)=\dot
 H_q^{2a+t}(\crnp)$, so  the estimate follows directly from Lemma 6.10.

Now let $t\in
 (-a+\frac1{q},1-a+\frac1q)$. When $u\in H_q^{a(2a+t)}(\crnp)$, it is
decomposed  into $u=w+z$ as in 
Corollary 3.3  (the $s$ there plays the role of $t$ here); then 
$$
w\in \dot H_q^{2a+t}(\crnp),\quad z=K^a_{0}\varphi \text{ with
}\varphi =\gamma _0^au\in B_q^{a+t-1/q}(\R^{n-1});  
$$
the norms being bounded in terms of
$\|u\|_{H_q^{a(2a+t)}(\crnp)}$. Let $\psi \in C_0^\infty (B_3)$ with
$\psi =1$ on $B_2$; then
$$
u=\psi u=w'+z', \text{ where } w'=\psi w\text{ and } z'=\psi z\text{ are
supported in }B_3.\tag6.34
$$

For $w'$ we have by Lemma 6.10:
$$
\aligned
\|x_n^{1-a+s}\nabla (w'/x_n^a)\|_{\dot H_q^{s+t}(\crnp)^n}
&\le C_1\|w'\|_{\dot H_q^{2a+t}(\crnp)}\\
&\le C_2\|w\|_{\dot H_q^{2a+t}(\crnp)}
\le C_3\|u\|_{ H_q^{a(2a+t)}(\crnp)}.
\endaligned 
$$
For $z'$ we shall use the  knowledge from Theorem 6.9 that $z$ satisfies
$$
\|x_n^{1-a+s}\nabla (z/x_n^a)\|_{\dot H_q^{s+t}(\crnp)^n}\le C_4
\|\varphi
\|_{B_q^{a+t-1/q}(\R^{n-1})}\le C_5 \|u\|_{ H_q^{a(2a+t)}(\crnp)}.\tag 6.35
$$
To estimate $z'$, write
$$
\aligned
\|&x_n^{1-a+s}\nabla (z'/x_n^a)\|_{\dot H_q^{s+t}(\crnp)^n}=
\|x_n^{1-a+s}\nabla (\psi z/x_n^a)\|_{\dot H_q^{s+t}(\crnp)^n
}\\
&\le \|x_n^{1-a+s}\psi \nabla( z/x_n^a)\|_{\dot
H_q^{s+t}(\crnp)^n}+\|x_n^{1-a+s}(\nabla \psi) z/x_n^a\|_{\dot
H_q^{s+t}(\crnp)^n}= I + II.
\endaligned
$$
Here $I\le C \|x_n^{1-a+s}\nabla( z/x_n^a)\|_{\dot
H_q^{s+t}(\crnp)^n}$, which is estimated in (6.35). For $II$ we have: 
$$
II\le C_1\|x_n^{1-a+s} z/x_n^a\|_{\dot
H_q^{s+t}(\crnp)}= C_2\|x_n^{1-a+s}K_0\varphi \|_{\dot
H_q^{s+t}(\crnp)}.
$$

The operator $K'=x_n^{1-a+s}K_0 $  acts as
$$
K'=c\Xi _+^{a-1-s}K_0,
$$
cf.\ (6.25).
Since $\varphi \in B_q^{a+t-1/q}(\R^{n-1})$, it is a fortiori in $
 B_q^{a+t-1-1/q}(\R^{n-1}) $ which is mapped into
$\ol H_q^{a+t-1}(\rnp)$
by
$K_0$. Here  $\ol H_q^{a+t-1}(\rnp)$
identifies with $\dot H_q^{a+t-1}(\rnp)$ since $a+t-1\in
(-\frac1{q'},\frac1q)$, and this space is mapped  by
$\Xi _+^{a-1-s}$ into $\dot H_q^{s+t}(\crnp)$. Altogether, we see that
$$
II\le C\|K'\varphi \|_{\dot H^{s+t}_q(\crnp)}\le
C'\|\varphi \|_{ B_q^{a+s+t-1/q}(\R^{n-1})};
$$
together with the estimate of $I$, this shows
$$
\|x_n^{1-a+s}\nabla (z'/x_n^a)\|_{\dot H_q^{s+t}(\crnp)^n}
\le C \|u\|_{ H_q^{a(2a+t)}(\crnp)}.
$$
Then (6.33) follows by adding the contributions from $w'$ and
$z'$.

We now have the estimate for $t\in
 (-a-\frac1{q'},1-a+\frac1q)$ except at the point $t=-a+\frac1q$; this
 is included by interpolation.
\qed 

\enddemo

Next, the estimates are established on sets with curved boundary by use of local coordinates:

\proclaim{Theorem 6.12}  Let $a\in (\frac12,1)$, $s\in [0,2a-1)$, $t\in
 (-a-\frac1{q'},1-a+\frac1q)$, $t\ge -s$,  $\tau \ge 1$,
 and let $\Omega $ be 
 a bounded
 $C^{1+\tau } $-domain. Assume that $s$ and $s+t$ are $<\tau -2a$. Then
$$
\|d^{1-a+s}\nabla (u/d^a)\|_{\dot H_q^{s+t} (\comega )}\le C
\|u\|_{H_q^{a(2a+t)}(\comega)}
.\tag6.36
$$

\endproclaim

\demo{Proof}
We use the definition of the spaces by
localization in the way recalled in Proposition 4.3; then the
estimates 
(6.33)
carry over to the spaces over $\Omega $.

As earlier noted, the differentiations $\partial_j$ are  by the
diffeomorphisms turned into first-order differentiations with $C^\tau
$-coefficients which preserve the estimates.
\qed 
\enddemo

An immediate application is to show gradient estimates of $u$ in
 terms of $f$, when $u$ is a solution of the homogeneous Dirichlet
 problem (1.1). 

\proclaim{Corollary 6.13}  Let $a\in(\frac12,1)$, $q\in (1,\infty )$,
 $\tau \ge 1$, let $\Omega $ be bounded and $C^{1+\tau }$, 
 and let $P$ be a strongly elliptic classical
 $2a$-order pseudodifferential operator satisfying Hypothesis {\rm
 2.1}. Let $s\in [0,2a-1)$, $t\in [0,1-a+\frac1q)$,  and assume that $s+t<\tau -2a$.
Let $u\in \dot H_q^a(\comega)$ be a solution of the Dirichlet problem
 {\rm (1.1)}. Then
$u$ and $f$ satisfy:
$$
\|d^{1-a+s}\nabla (u/d^a)\|_{\dot H_q^{s+t} (\comega )}\le C
(\|f\|_{\ol H_q^t(\Omega )}+\|u\|_{L_q(\Omega )}).
\tag6.37
$$
When there is uniqueness of solution, the term
 $\|u\|_{L_q(\Omega )}$ can be omitted.
\endproclaim

\demo{Proof} Estimates with $\|u\|_{ H_q^{a(2a+t)}(\comega)}$ in the
right-hand side have been shown in Theorem 6.12. We
have from Theorem 2.2, cf.\ (2.10), that for solutions in the case with
unique solvability, $\|u\|_{
H_q^{a(2a+t)}(\comega)}$ in the right-hand side can be replaced by
$\|f\|_{\ol H_q^{t}(\Omega)}$.
In the general case we can apply the right-hand side inequality in
(2.11). \qed

\enddemo

\subsubhead 6.4 Estimates of $\nabla( u/d^a)$
 in H\"older-Zygmund spaces \endsubsubhead

The gradient estimates in H\"older-Zygmund spaces will be shown by
following the pattern  set up for Bessel-potential spaces.

\proclaim{Lemma 6.14} Let $a\in (\frac12,1)$, $s \in [0,2a-1)$, $t\ge 
 -s$.

$1^\circ$ When $w\in \dot C_*^{2a+t}(\crnp)$ and is supported in $\ol B $,
$$
\|x_n^{1-a+s }\nabla (w/x_n^a)\|_{\dot C_*^{s+t} (\crnp)^n}\le
C\|w\|_{\dot C_*^{2a+t}(\crnp)}. 
\tag6.38
$$

$2^\circ$ Consider a bounded $C^{1+\tau }$-domain with $\tau \ge
1$, and assume that $s+t$ and $t$ are $<\tau -2a$. Here  $w\in \dot C_*^{2a+t}(\comega)$ satisfies
$$
\|d^{1-a+s }\nabla (w/d^a)\|_{\dot C_*^{s+t} (\comega)^n}\le
C\|w\|_{\dot C_*^{2a+t}(\comega)}. 
\tag6.39
$$
\endproclaim

\demo{Proof}
 $1^\circ$. Since $w\in \dot C_*^{2a+t}(\crnp)$, $\nabla w\in \dot
C_*^{2a-1+t}(\crnp)^n$, and hence by  (5.24),
$$
\nabla w\in \dot C_*^{2a-1+t}(\crnp)^n\subset x_n^{2a-1-s }\dot
C_*^{s +t}(\crnp)^n,  
$$
with norm dominated by that of $w$ in $\dot C_*^{2a+t}(\crnp)$.
Then for $w/x_n^a$, 
$$
x_n^{1-a+s }\nabla (w/x_n^{a})=
x_n^{1-2a+s }\nabla w -ax_n^{-2a+s }(0,\dots,0,w)
\in \dot C_*^{s+t} (\crnp)^n
$$
by (5.24), implying 
(6.38).

$2^\circ$. The estimates carry over to curved sets by use of local
coordinates as in Proposition 4.3.
\qed
\enddemo

This lemma will take care of the part in $\dot C_*^{2a+t}(\comega)$ of a
function in $C_*^{a(2a+t)}(\comega)$.

For functions $z=K^a_0\varphi $, we have found the structure of the Poisson-like
operators arising from gradient formulas in Theorem 6.9; now we just have
to study their mapping properties in $C_*^t$-spaces.

\proclaim{Theorem 6.15} Let $a\in (0,1)$. For $z=K_0^a\varphi $
on $\rnp$, there holds when
$0< a+t<1$, $s\ge 0$,
$$
\|x_n^{1-a+s}\nabla (z/x_n^a)\|_{\dot C_*^{s+t}(\crnp)^n}\le  C\|\varphi
\|_{C_*^{a+t}(\R^{n-1})}.\tag 6.40
$$
\endproclaim

\demo{Proof} Recall from (6.27) that $x_n^{1-a+s}\nabla (z/x_n^a)=K^{(n)}\varphi $
with $K^{(n)}$ written as in (6.28). Note that $\ang{D'}$ maps  $\varphi
\in C_*^{a+t}(\R^{n-1})$ to $C_*^{a+t-1}(\R^{n-1})$, which is mapped to
$\ol C_*^{a+t-1}(\rnp)$ by $K_0$.  Here $\ol C_*^{a+t-1}(\rnp)$
identifies with $\dot C_*^{a+t-1}(\rnp)$ when $a+t\in (0,1)$ (cf.\ the
lines after (2.19)).
Applying $\Xi _+^{a-1-s}$, we get $K^{(n)}\varphi \in \dot C_*^{s+t}(\crnp)$.
Altogether, 
$K^{(n)}$ maps
 $$
K^{(n)}\colon C_*^{a+t}(\R^{n-1})\to \dot C_*^{s+t}(\crnp),
\text{ when $a+t\in (0,1)$, $s\ge 0$} .\tag6.41
$$

The operators $K^{(j)}$ for $j<n$ are treated in a similar way,
resulting in (6.40). \qed
\enddemo

This is combined with the estimates in Lemma 6.14:

\proclaim{Theorem 6.16}  Let $a\in (\frac12,1)$, $s\in [0,2a-1)$, $t\in
 (-a,1-a)$ and $t\ge -s$.
For $u\in C_*^{a(2a+t)}(\crnp) $ with $\supp u\subset \ol B$, there holds:
$$
\|x_n^{1-a+s}\nabla (u/x_n^a)\|_{\dot C_*^{s+t} (\crnp )^n}\le C \|u\|_{C_*^{a(2a+t)}(\crnp)}.\tag6.42
$$

\endproclaim

\demo{Proof} The proof is completey analogous to the proof of Theorem 6.11, so we
leave out the details. \qed

\enddemo

As in the case of $H^t_q$-spaces, the estimates carry over to bounded
$C^{1+\tau }$-domains by localization as described in Proposition 4.3:

\proclaim{Theorem 6.17}  Let $a\in (\frac12,1)$, $s\in [0,2a-1)$, $t\in
 (-a,1-a)$, $t\ge -s$,  $\tau \ge 1$,
 and let $\Omega $ be 
 a bounded
 $C^{1+\tau } $-domain. Assume that $s$ and $s+t$ are $<\tau -2a$. Then
$$
\|d^{1-a+s}\nabla (u/d^a)\|_{\dot C_*^{s+t} (\comega )^n}\le C \|u\|_{C_*^{a(2a+t)}(\comega)}.\tag6.43
$$
\endproclaim

Finally, we address the question of  gradient estimates for
 solutions of (1.1) in H\"older
spaces.
In the smooth case, the solvability result
from \cite{G14}
recalled in
Theorem 2.3 $1^\circ$ allows the conclusion:

\proclaim{Theorem 6.18}  Let $a\in (\frac12,1)$, let $\Omega $ be bounded
 $C^\infty $ and let $P$ satisfy Hypothesis {\rm 2.1} with $\tau =\infty
 $. Let $s\in [0,2a-1)$, $t\in
 (-a,1-a)$ and $t\ge -s$. Let $u$
 solve {\rm (1.1)}.
Then
$$
\aligned
\|d^{1-a+s}\nabla (u/d^a)\|_{\dot C_*^{s+t} (\comega )^n}&\le C
(\|f\|_{\ol C_*^t(\Omega)}+\|u\|_{L_\infty (\Omega )}),\\
\|d^{1-a+s}\nabla (u/d^a)\|_{\dot C_*^s (\comega )^n}&\le C
(\|f\|_{L_\infty (\Omega)}+\|u\|_{L_\infty (\Omega )}).
\endaligned\tag6.44
$$
When there is uniqueness of solution, the terms  $\|u\|_{L_\infty
(\Omega)}$ can be omitted.
\endproclaim

\demo{Proof} We know that the solution $u$ is in
$C_*^{a(2a+t)}(\comega)$, and (6.43) holds.
The replacement of the  right-hand side in (6.43) by the expressions with
$f$ goes as accounted for in Corollary 6.3.
\qed

\enddemo

In nonsmooth cases, we shall derive
the gradient estimates in $C_*^s$-spaces for solutions to (1.1) simply by applying
Sobolev embedding and other embedding properties to the
$H_q^s$-results in Corollaries 6.7 and 6.13. It is well-known that there
holds for bounded $C^\sigma $-domains $\Omega $:
$$
\ol C_*^s(\Omega )\subset \ol H_q^{s-\varepsilon /2}(\Omega )\subset
\ol C_*^{s-\varepsilon -n/q}(\Omega ),\text{ when }0\le s-\varepsilon
-n/q<s<\sigma .\tag 6.45
$$
Then we finally obtain:

\proclaim{Theorem 6.19}  Let $a\in(\frac12,1)$, 
 $\tau \ge 1$, let $\Omega $ be bounded and $C^{1+\tau }$, 
 and let $P$ be a strongly elliptic classical
 $2a$-order pseudodifferential operator satisfying Hypothesis {\rm
 2.1}. Let $s\in [0,2a-1)$, $t>0$,
and assume that $s+t<\tau -2a$.
Let $u\in \dot C_*^a(\comega)$  be a solution of the Dirichlet problem {\rm (1.1)}. Then
$u$ and $f$ satisfy for small $\varepsilon >0$:
$$
\aligned
\|d^{1-a}\nabla u\|_{\ol C_*^{s+t-\varepsilon } (\Omega  )^n}&\le C
(\|f\|_{\ol C_*^t(\Omega )}+\|u\|_{L_\infty (\Omega )}),\\
\|d^{1-a}\nabla u\|_{\ol C_*^{2a-1-\varepsilon } (\Omega  )^n}&\le C
(\|f\|_{L_\infty (\Omega )}+\|u\|_{L_\infty (\Omega )}).
\endaligned
\tag6.46
$$
When moreover  $t\in (0,1-a)$, $s>0$,
$$
\aligned
\|d^{1-a+s}\nabla (u/d^a)\|_{\dot C_*^{s+t-\varepsilon } (\comega )}&\le C
(\|f\|_{\ol C_*^t(\Omega )}+\|u\|_{L_\infty (\Omega )}),\\
\|d^{1-a+s}\nabla (u/d^a)\|_{\dot C_*^{s-\varepsilon } (\comega )}&\le C
(\|f\|_{L_\infty (\Omega )}+\|u\|_{L_\infty (\Omega )}).
\endaligned
\tag6.47
$$
\endproclaim

\demo{Proof} Corollary 6.7 implies the inequality for $s\le 2a-1$:
$$
\|d^{1-a}\nabla u\|_{\ol H_q^{s+t} (\Omega  )^n}\le C
(\|f\|_{\ol H_q^t(\Omega )}+\|u\|_{L_q(\Omega )}).
$$
By using (6.45) we find:
$$
\aligned
\|d^{1-a}\nabla u\|_{\ol C_*^{s+t-\varepsilon } (\Omega  )^n}&\le C_1
\|d^{1-a}\nabla u\|_{\ol H_q^{s+t-\varepsilon/2 } (\Omega  )^n}
\le
C_2(\|f\|_{\ol H_q^{t-\varepsilon /2}(\Omega )}+\|u\|_{L_q(\Omega
)})\\
&\le C_3(\|f\|_{\ol C_*^{t-\varepsilon /4+n/q}(\Omega )}+\|u\|_{L_q(\Omega
)}) \le C_4(\|f\|_{\ol C_*^{t}(\Omega )}+\|u\|_{L_\infty (\Omega
)}),
\endaligned
$$
where $q$ is taken so large that $n/q\le \varepsilon /4$. This shows
the first estimate in (6.46). The second estimate follows, taking
$s=2a-1-\varepsilon/2 $, $t=0$, from
$$
\aligned
\|d^{1-a}\nabla u\|_{\ol C_*^{2a-1-\varepsilon } (\Omega  )^n}&\le C_1
\|d^{1-a}\nabla u\|_{\ol H_q^{2a-1-\varepsilon/2 } (\Omega  )^n}\\
&\le C_2(\|f\|_{L_q(\Omega )}+\|u\|_{L_q(\Omega
)}) \le C_3(\|f\|_{L_\infty (\Omega )}+\|u\|_{L_\infty (\Omega
)}).
\endaligned
$$

From Corollary 6.13 we have when $0\le t<1-a+\frac1q$, $0\le s<2a-1 $,
$$
\|d^{1-a+s}\nabla (u/d^a)\|_{\dot H_q^{s+t} (\comega )}\le C
(\|f\|_{\ol H_q^t(\Omega )}+\|u\|_{L_q(\Omega )}).
$$
Here an application of (6.45) gives when $t>0$, $s>0$:
$$
\aligned
\|d^{1-a+s}\nabla (u/d^a)\|_{\dot C_*^{s+t-\varepsilon } (\Omega  )^n}&\le C_1
\|d^{1-a+s}\nabla (u/d^a)\|_{\dot H_q^{s+t-\varepsilon/2 } (\Omega  )^n}\\
&\le
C_2(\|f\|_{\ol H_q^{t-\varepsilon /2}(\Omega )}+\|u\|_{L_q(\Omega
)})\\
&\le C_3(\|f\|_{\ol C_*^{t-\varepsilon /4+n/q}(\Omega )}+\|u\|_{L_q(\Omega
)})\\
&\le C_4(\|f\|_{\ol C_*^{t}(\Omega )}+\|u\|_{L_\infty (\Omega
)}),
\endaligned
$$
for $n/q\le \varepsilon /4$, showing the first estimate in (6.47). For the second estimate, we take
$t=0$ and observe that for $s>0$,
$$
\aligned
\|d^{1-a+s}\nabla (u/d^a)\|_{\dot C_*^{s-\varepsilon } (\Omega  )^n}&\le C_1
\|d^{1-a+s}\nabla (u/d^a)\|_{\dot H_q^{s } (\Omega  )^n}\\
&\le
C_2(\|f\|_{L_q(\Omega )}+\|u\|_{L_q(\Omega
)})\\
&\le C_3(\|f\|_{L_\infty (\Omega )}+\|u\|_{L_\infty (\Omega
)}) . \qed
\endaligned 
$$

\enddemo

The second line in (6.47) implies in particular, since $\dot
C_*^\varepsilon (\comega)\subset L_\infty (\Omega )$,
$$
\|d^{1-a+\delta }\nabla (u/d^a)\|_{L_\infty (\Omega ) )}\le C
(\|f\|_{L_\infty (\Omega )}+\|u\|_{L_\infty (\Omega )^n})\text{ for }\delta >0,\tag6.48
$$
extending (6.2) to hold for general operators $P$.

\subhead Appendix A. Localization in H\"older-Zygmund spaces \endsubhead

It was shown in \cite{AG23, Prop.\ 4.5} that the spaces  
 $H_q^{\mu (s)}(\crnp)$  are invariant under multiplication by a 
 $C^\sigma $-function $\varphi $ under suitable conditions on $\sigma
 $ relative to  $\mu $ and $s$:
 $$\aligned
u\in H_q^{\mu (s)}(\crnp)
&\implies \varphi u\in H_q^{\mu (s)}(\crnp)\text{ when}\\
&\text{either }\mu \ge 0, \quad \sigma >\mu,  \quad \mu
-\tfrac1{q'}<s<\sigma ,\\
&\text {or }\mu \in (-1,0), \quad \sigma >\mu +1,  \quad \mu -\tfrac1{q'}<s<\sigma -1.\
\endaligned
\tag A.1$$

It was used in the subsequent Theorem 4.6 that a related statement
holds for $C^s_*$-spaces, but details were not given. We shall provide
the details here.

Some auxiliary observations:

Recall that $B^s_{\infty
,\infty }(\rn)=C_*^s(\rn)$. The following embedding properties
hold for $s\in\R$, $\varepsilon >0$, when $\psi \in C_0^\infty (\rn)$:
$$
H_q^s(\rn)\subset C_*^{s-\frac nq -\varepsilon
}(\rn), \quad 
\psi  C_*^{s}(\rn )\subset  H_q^{s-\varepsilon }(\rn ).
\tag A.2
$$
Based on (A.2) there is an unsharp result for general parameters:

\proclaim{Lemma A.1} Let $\mu ,\sigma ,s$ satisfy 
either $\mu \ge 0, \sigma >\mu, s\in ( \mu
-1,\sigma )$, or $\mu \in (-1,0), \sigma >\mu +1, s\in ( \mu
-1,\sigma -1)$,
and let $\varphi \in
C^\sigma (\rn)$ with compact support. For a large ball $B_R=\{|x|<R\}$  containing $\supp
\varphi $, let  $\psi \in
C_0^\infty (B_R)$ with $\psi =1$ on a neighborhood of $\supp \varphi $.
Then for $s-\mu >\varepsilon >0$,
$$
\varphi C_*^{\mu (s)}(\crnp) \subset \psi  C_*^{\mu (s-\varepsilon
)}(\crnp)\subset  C_*^{\mu (s-\varepsilon
)}(\crnp).\tag A.3
$$
\endproclaim

\demo{Proof}  
Let $u\in C_*^{\mu (s)}(\crnp)$, then $\varphi u=\varphi
\psi u$. Now
$$
\varphi \psi u=\varphi u\in \varphi C_*^{\mu (s)}(\crnp)=\varphi \Lambda _+^{-\mu }e^+\ol
C_*^{s-\mu }(\rnp),
$$
so there is a $v\in \ol
C_*^{s-\mu }(\rnp)$ such that $\varphi u= \varphi \Lambda _+^{-\mu }e^+v$. Moreover,
$$
\varphi u= \varphi  \Lambda _+^{-\mu }e^+v=\varphi  \Lambda _+^{-\mu
}\psi e^+v+\varphi  \Lambda _+^{-\mu }(1-\psi )e^+v.\tag A.4
$$
For the first term we note that
$$
\aligned
\psi e^+v&\in \psi e^+\ol
C_*^{s-\mu }(\rnp)\subset e^+\{w\in \ol
C_*^{s-\mu }(\rnp)\mid \supp w\subset B_R)\\
&\subset e^+\{w\in \ol
H_q^{s-\mu -\varepsilon }(\rnp)\mid \supp w\subset B_R),
\endaligned
$$
in view of (A.2). Then the localization rule for $H_q^s$-spaces (A.1)
(applied for a large enough $q$) shows that
$$
\aligned
\varphi  \Lambda _+^{-\mu }\psi e^+v
&\in  \varphi  \Lambda _+^{-\mu } e^+H_q^{s-\mu -\varepsilon
}(\rnp)\\
&=\varphi H_q^{\mu (s-\varepsilon )}(\crnp)\subset  H_q^{\mu
(s-\varepsilon )}(\crnp)\subset C_*^{\mu (s-\varepsilon -\frac nq-\varepsilon ')}(\crnp);
\endaligned
$$
in the last step we used (A.2).

For the second term in (A.4) we note that $\varphi  \Lambda _+^{-\mu
}(1-\psi )$ is of order $-\infty $, since $1-\psi $
vanishes on a neighborhood of $\supp \varphi $, so the term is in $C^\sigma  (\rn)$. It is
supported in $\crnp \cap B_R$, since the other terms are so; in
particular, it 
is in $\dot C_*^\sigma (\crnp)\subset C_*^{\mu (\sigma )}(\crnp)$.

Since we can take $q$ arbitrarily large,  (A.3) follows. \qed
\enddemo

This shows a suboptimal form (because of the loss of $\varepsilon $)
of the localization property we need. We shall now show that $\varepsilon $ can be removed under suitable
hypotheses on the parameters.

Recall the  well-known property: For $\varphi \in C^\sigma (\rn)$,
$$
\varphi C_*^t(\rn)\subset C_*^t(\rn)\text{ when }|t|<\sigma .\tag A.5
$$

It will be useful to have the following composition rule from \cite{A05, Th.\ 3.6}, applied
to $\Lambda _+^t  $ and the multiplication by a function $\varphi
\in C^\sigma (\rn)$:

Let $\sigma >0$ and  $\theta \in (0,\sigma )$
($[\theta ]=$ largest integer $<\theta $), $t \in \R$ and  $\varphi \in C^\sigma (\rn)$. Then 
$$
\Lambda _+^{t }\varphi -\sum_{0\le |\alpha |\le [\theta ]}\OP(\tfrac1{\alpha
!}\partial_\xi ^\alpha \lambda _+^{t }(\xi )D_x^\alpha \varphi
(x))=R_{[\theta ]}\colon C_*^{s+t -\theta }(\rn)\to C_*^{s}(\rn)\tag A.6
$$
is bounded when $|s|<\sigma $, $s-\theta  >-\sigma  $, $-\sigma
+\theta <s+t <\sigma $.

These parameter conditions simplify to 
$$
\aligned
-\sigma +\theta &<s<\sigma , \text{ when }t  \le 0,\\
-\sigma +\theta &<s<\sigma -t , \text{ when }t\ge 0.
\endaligned
\tag A.7
$$

\proclaim{Theorem A.2} Let $\mu >-1$, $\sigma >0$, $\varphi \in C^\sigma (\rn)$. The localization
 property
 $$
u\in C_*^{\mu (s)}(\crnp)
\implies \varphi u\in C_*^{\mu (s)}(\crnp)\tag A.8
$$ holds under the
following conditions:

$1^\circ$ $s\in (\mu -1,\mu )$ and $\sigma \ge \max\{\mu ,1-\mu\}$.

$2^\circ$ $s \in [\mu ,\mu +1)$ and $s<\sigma $, where $\sigma \ge\mu $
if $\mu \ge 0$,  $\sigma \ge 1-\mu $
if $\mu < 0$.

$3^\circ$ $\supp\varphi $ is compact, $s\ge \mu +1$, and one of the conditions {\rm (A.9)-(A.11)} holds: 
$$
\align
&\mu \ge 1,\quad \sigma >\mu  , \quad  s<\sigma -1, \tag A.9\\
& \mu \in [0,1),\quad \sigma > \mu +1, \quad  s <\sigma
-2.\tag A.10\\
&\mu \in (-1,0), \quad \sigma >\mu +3,\quad s < \sigma
 -3.\tag A.11
\endalign
$$

\endproclaim

\demo{Proof}
$1^\circ$. 
For $s\in (\mu
-1,\mu )$, $$
C_*^{\mu (s)}(\crnp)=\Lambda _+^{-\mu }e^+\ol C_*^{s-\mu }(\rnp)
=\Lambda _+^{-\mu }\dot  C_*^{s-\mu }(\crnp)=
\dot C_*^s(\crnp).\tag A.12
$$
Then (A.8) holds in view of (A.5), since the conditions assure that $s\in (\mu -1,\mu )$ implies $|s|<\sigma $. 

$2^\circ$. Let $s\in (\mu , \mu +1)$.
Note that by definition, there holds in general for $s-\mu >-1$:
$$
u\in C_*^{\mu (s)}(\crnp)\iff u=\Lambda _+^{-\mu }v,\text{ where }v\in
e^+\ol C_*^{s-\mu }(\rnp ).\tag A.13
$$
Then
$$
\varphi u= \varphi \Lambda _+^{-\mu }v=\Lambda _+^{-\mu}(\varphi v)+[\varphi ,\Lambda _+^{-\mu } ]v.
\tag A.14
$$
Here  $\varphi v\in e^+\ol C_*^{s-\mu }(\rnp )$ when $|s-\mu
|<\sigma$ by the rule (A.5); then the first term in the last entry of
(A.14) is in $C_*^{\mu (s)}(\crnp)$.

For the second term, we use the composition formula (A.6) with $t=-\mu $,
$\theta \in (0,1)$. Here $[\theta ]=0$ and
$$
R_0=\Lambda _+^{-\mu }\varphi -\OP(\varphi (x)\lambda _+^{-\mu }(\xi
))=-[\varphi ,\Lambda_+^{-\mu }]\colon C_*^{s-\mu -\theta }(\rn) \to
C_*^{s}(\rn) 
$$
when (A.7) holds with $t=-\mu $.
For $s\in (\mu ,\mu +1)$, we can write $s=\mu
+\theta -\varepsilon $ for a $\theta \in (0,1)$ and a small
$\varepsilon $.
Then
$$
R_0\colon C_*^{-\varepsilon  }(\rn)\to C_*^{s}(\rn),
$$
and since $R_0$ preserves support in $\crnp$,
$$
R_0\colon \dot C_*^{-\varepsilon  }(\crnp)\to \dot C_*^{s}(\crnp)\subset
C_*^{\mu (s)}(\crnp).
$$
Now $v\in e^+\ol C_*^{s-\mu }(\rnp)\subset e^+\ol C_*^{s-\mu -\theta
}(\rnp)=e^+\ol C_*^{\,-\varepsilon  }(\rnp)= \dot C_*^{-\varepsilon
}(\crnp)$, so it is mapped by $R_0$ into $C_*^{\mu (s)}(\crnp)$, as
desired.

The value $s=\mu $ is included by interpolation with the case
$1^\circ$.
This ends the proof of $2^\circ$.

The method can be pursued for higher values of $s$ and $\theta $, but
does not seem to give a better result for $s>\mu +1$ than the
method we shall use to show $3^\circ$.

$3^\circ$. The statements here  allow arbitrarily large
$s$, dominated by $\sigma $ in a linear way. 

{\it The method for }(A.9).
 For a large ball $B_R=\{|x|<R\}$  containing $\supp
\varphi $, let  $\psi ,\psi _1 \in
C_0^\infty (B_R)$ with $\psi _1=1$ on a neighborhood of $\supp \psi
 $ and  $\psi =1$ on a neighborhood of $\supp \varphi $.

By (A.13),
$u\in C_*^{\mu (s)}(\crnp)$ means that $\Lambda
_+^\mu u \in e^+\ol C_*^{s-\mu }(\rnp)$. 
So we must show that $\Lambda _+^\mu (\varphi u)\in e^+\ol
C_*^{s-\mu }(\rnp)$ (or, equivalently, $r^+\Lambda _+^\mu (\varphi u)\in \ol
C_*^{s-\mu }(\rnp)$). Here $$
\Lambda _+^\mu (\varphi u)=\varphi \Lambda _+^\mu
u +[\Lambda _+^\mu ,\varphi ]u,\tag A.15
$$ where  $\varphi \Lambda _+^\mu
u  \in e^+\ol C_*^{s-\mu }(\rnp)$ in view of (A.5), when $s-\mu
<\sigma $, i.e.\ $s<\sigma +\mu $. It remains to treat the term
$[\Lambda _+^\mu ,\varphi ]u$.

Since $\varphi u=\varphi \psi u$, we can replace $u$ by $\psi u$ throughout
the formula (A.15) (or
assume from the start that $u$ is supported in $B_R$).

The
commutator $[\Lambda _+^\mu ,\varphi ]$ was studied in \cite{AG23,
p.1317}, and we have, continuing the calculations from there:
$$
\aligned
[\Lambda _+^\mu ,\varphi ]&=\OP(a(x,y,\xi )),\text{ where }\\
a(x,y,\xi )&=\lambda
_+^\mu (\xi )(\varphi (x)-\varphi (y))=\lambda _+^\mu (\xi
){\sum}_{j\le n}(x_j-y_j)\varphi _j(x,y)
\endaligned\tag A.16
$$
with $\varphi _j\in C^{\sigma -1}(\R^{2n})$. Moreover, by the Fourier
integral definitions, $\OP(a(x,y,\xi
))=\OP(a'(x,y,\xi ))$, where
$$
a'(x,y,\xi )={\sum}_{j\le n}D_{\xi _j}\lambda _+^\mu (\xi )\varphi
_j(x,y)\in C^{\sigma -1}S^{\mu -1}(\R^{2n}\times \rn).\tag A.17
$$
Denote the operator $\OP(a'(x,y,\xi ))=A$.

We want to apply \cite{AG23, Th.\ 5.13} to the mapping $A$ in $H_q^s$-spaces, and use
(A.2) in a similar way as in Lemma A.1 to reach $C_*^s$-spaces.
First we note that the order $\mu -1$ is lower than the parameter $\mu $, not
allowed in Th.\ 5.13, so a
modification is needed. Here we shall use that 
$$
H_q^{\mu (s)}(\crnp)\subset  H_q^{(\mu -1)(s)}(\crnp), 
$$
so that we just have to handle a mapping from $ H_q^{(\mu -1)(s)}(\crnp)$.

Assume $\mu \ge 1$. By \cite{AG23, Th.\ 5.13}, 
$$
r^+A\colon  H_q^{(\mu -1)(s'+\mu -1)}(\crnp)\to \ol H_q^{s'}(\rnp),  
$$
when $\sigma -1>0$, $\mu -1<\sigma -1$, $0\le s'<\sigma -1-(\mu -1)$, i.e.,
$
\sigma >1,\mu <\sigma ,0\le s'<\sigma -\mu .
$
With $s=s'+\mu -1$, this means that
$$
r^+A\colon  H_q^{(\mu -1)(s)}(\crnp)\to \ol H_q^{s-\mu
+1}(\rnp), \tag A.19 
$$
when 
$$
\sigma > \mu  ,\quad\mu -1\le s<\sigma -1.\tag A.20
$$

With $\psi $, $\psi _1$ defined above, $u$  can be replaced by $\psi u$, $\supp \psi
\subset B_R$, and moreover,
$$
\aligned
\psi u\in\psi C_*^{(\mu -1)(s)}(\crnp)&=\psi \Lambda _+^{\mu -1}e^+\ol
C_*^{s-\mu +1}(\rnp)\\
&\subset \psi \Lambda _+^{\mu -1}\psi _1e^+\ol
C_*^{s-\mu +1}(\rnp)+ \psi \Lambda _+^{\mu -1}(1-\psi _1)e^+\ol
C_*^{s-\mu +1}(\rnp).
\endaligned
$$
Here since $\psi $ and $1-\psi _1$ have disjoint supports,  $\psi \Lambda _+^{\mu -1}(1-\psi _1)$ is a smoothing operator
mapping everything into $C^\infty (\rn)$, in the present case to
functions supported in $\crnp\cap B_R$, lying in $\psi _1 \dot C^\infty
(\crnp)$.
Now by (A.2),
$$
\aligned
\psi u&\in \psi \Lambda _+^{\mu -1}\psi _1e^+\ol
C_*^{s-\mu +1}(\rnp)+\psi _1 \dot C^\infty
(\crnp)\\
&\subset \psi \Lambda _+^{\mu -1}e^+\ol
H_q^{s-\mu +1-\varepsilon }(\rnp)+\psi _1 \dot C^\infty
(\crnp)=\psi H_q^{(\mu
-1)(s-\varepsilon )}(\crnp)+\psi _1 \dot C^\infty
(\crnp)\\
&\subset  H_q^{(\mu
-1)(s-\varepsilon )}(\crnp).
\endaligned
$$
Here $r^+A$ can be applied according to (A.19), mapping into $$
r^+A\psi u\in \ol
H_q^{s-\mu +1-\varepsilon}(\rnp)\subset \ol C_*^{s-\mu +1-\varepsilon
-\frac nq-\varepsilon '}(\rnp)\subset 
\ol C_*^{s-\mu}(\rnp),
$$
when $q$ is so large and $\varepsilon ,\varepsilon '$ are so small
that $\varepsilon +\varepsilon '+\frac nq\le 1$. We conclude that
$$
r^+A\colon  \psi C_*^{\mu (s)}(\crnp) \subset \psi C_*^{(\mu -1)(s)}(\crnp)
\to  \ol C_*^{s-\mu
}(\rnp),\tag A.21
$$
when (A.20) holds.

This accounts for the second term in (A.15), so we have altogether
found that $r^+\Lambda _+^\mu (\varphi u)\in \ol C_*^{s-\mu }(\rnp)$,
as was to be shown. Thus the multiplication by $\varphi $ maps
 $C_*^{\mu (s)}(\crnp) $ into itself.

{\it The method for }(A.10). Now let $\mu \in (0,1)$. With $\sigma '=\sigma -1>0$, $\mu '=\mu
-1\in (-1,0)$, we have that $a'(x,y,\xi )\in C^{\sigma '}S^{\mu '}(\R^{2n}\times
\rn)$. Here \cite{AG23, Cor.\ 5.14} shows that 
$$
r^+\OP(a')\colon  H_q^{\mu'(s'+\mu ')}(\crnp)\to \ol H_q^{s'}(\rnp),  
$$
when $\mu '<\sigma '-1$, $0\le s'<\sigma '-\mu '-1$. For
 $s=s'+\mu '$, the conditions mean in terms of $\mu $ and
$\sigma $  that
$$
\sigma -1 > \mu  ,\quad \mu -1\le s <\sigma -2;\tag A.22
$$
then (A.19) holds.
We can conclude as above by taking  $q$ sufficiently large that (A.21)
holds under the conditions (A.22); this shows the rest of the proof that (A.8)
holds under the conditions (A.22) (hence (A.9)) in the case $\mu \in (0,1)$.

{\it The method for }(A.11).
Here we can  use results from
the case
$\mu \in (0,1)$, as follows: Since
$$
\aligned
C_*^{\mu (s)}(\crnp)&=\Xi _+^{-\mu }e^+\ol C_*^{s-\mu }(\rnp)=\Xi
_+^1\Xi _+^{-\mu -1}e^+\ol C_*^{s-\mu }(\rnp)
=\Xi _+^1C_*^{(\mu +1)(s+1)}(\crnp)\\
&=(\ang{D'}+\partial_n)  C_*^{(\mu
+1)(s+1)}(\crnp)\subset C_*^{(\mu +1)(s)}(\crnp)+\partial_nC_*^{(\mu +1)(s+1)}(\crnp),
\endaligned\tag A.23
$$
a function $u\in  C_*^{\mu (s)}(\crnp)$ can be written as
$u=v+\partial_nw$, where 
$$
v\in C_*^{(\mu +1)(s)}(\crnp), \quad w\in C_*^{(\mu +1)(s+1)}(\crnp).
$$
For $v$, the rule (A.10)
gives
$$
\varphi v\in C_*^{(\mu +1)(s)}(\crnp)\text{
when }\sigma -1 >\mu +1,\; \mu \le s<\sigma -2.\tag A.24
$$
For $w$,
$$
\varphi \partial_nw=\partial_n(\varphi w)-\varphi _n w, \text{ where
}\varphi _n=\partial_n\varphi \in C^{\sigma -1}(\rn).
$$
Here $\partial_n(\varphi w)$ is treated by
use of the  rule (A.10) inside the bracket:
$$\aligned
\partial_n\big[\varphi C_*^{(\mu
+1)(s+1)}(\crnp)\big]&\subset \partial_n C_*^{(\mu
+1)(s+1)}(\crnp)\\
&\subset \Xi _+^1 C_*^{(\mu
+1)(s+1)}(\crnp)+\ang{D'} C_*^{(\mu +1)(s+1)}(\crnp)\\
&\subset  C_*^{\mu (s)}(\crnp)+   C_*^{(\mu +1)(s)}(\crnp)=C_*^{\mu
(s)}(\crnp),\\
&\text{ when }\sigma -1>\mu +1,\; \mu \le s+1<\sigma -2.
\endaligned \tag A.25
$$
For $\varphi _nw$,
we have by (A.10) with $\sigma '=\sigma -1$, $\mu '=\mu +1$,
$$
\aligned
\varphi _nC_*^{(\mu +1)(s+1)}(\crnp)&\subset
\varphi _nC_*^{(\mu +1)(s)}(\crnp)\subset
C_*^{(\mu
+1)(s)}(\crnp) \\
&\text{ when }\sigma '>\mu '+1,\; \mu '-1\le s<\sigma '
-2;
\endaligned
\tag A.26
$$
these conditions are in terms of the original values:
$$
\sigma >\mu +3,\quad \mu \le s<\sigma -3.\tag A.27
$$

The contributions in (A.24)-(A.26) show together that 
$$
\varphi C_*^{\mu (s)}(\crnp)\subset  C_*^{\mu (s)}(\crnp)
$$
when (A.27) holds, so (A.11) can be concluded.
\qed 
\enddemo

Admittedly, the conditions found in this theorem for high $s$ are not
very elegant; it is possible that they may be improved by other methods.

The main body of the present paper 
considers the case $\mu =a\in (0,1)$, so let
us list the consequences for this case:

\proclaim{Corollary A.3} Let $a\in (0,1)$, $\varphi \in C^\sigma
(\rn)$. There holds
$$
u\in C_*^{a (s)}(\crnp)\implies \varphi u\in C_*^{a (s)}(\crnp),\tag A.28
$$
when one of the following conditions is satisfied:
$$
\align
&a-1<s<a,\quad \sigma \ge \max\{a ,1-a\},\tag A.29\\
&a\le s< a+1,\quad \sigma \ge a ,\quad s<\sigma ,\tag A.30\\
&a +1\le s <\sigma -2, \text{ $\supp\varphi $ is compact.} 
\tag A.31
\endalign
$$
\endproclaim

\demo{Proof} Follows directly from $1^\circ$-$3^\circ$ in Theorem A.2.
\qed
\enddemo

Note that (A.30) implies the simpler statement that (A.28) holds when
$\sigma  \ge a+1$ and $s\in [a,a+1)$; but this is more demanding of
$\sigma $.

\enddocument